\documentclass[a4paper,onecolumn, superscriptaddress,10pt, accepted=2022-05-30, issue=2, volume=4, shorttitle=papers/compositionality-4-2]{compositionalityarticle}
\pdfoutput=1
\usepackage{amssymb}
\usepackage{amsmath}
\usepackage{array}
\usepackage{bbm}
\usepackage{enumerate}
\usepackage[latin1]{inputenc}
\usepackage[T1]{fontenc}
\usepackage{shadow}
\usepackage[all,cmtip,2cell]{xy}
\UseAllTwocells
\usepackage{graphicx}
\usepackage{graphbox}
\usepackage{caption}
\usepackage{subcaption}
\usepackage{float}
\usepackage{amsthm}
\usepackage{multicol}
\usepackage[]{hyperref}
\usepackage{tikz}
\usepackage{tikz-cd}
\usepackage{longtable}
\usepackage{color}
\usetikzlibrary{matrix,arrows,decorations.pathmorphing}
\usepackage{pdfpages}

\newtheorem{thm}{Theorem}[section]
\newtheorem{df}[thm]{Definition}
\newtheorem{p}[thm]{Proposition}
\newtheorem{lmm}[thm]{Lemma}

\newtheorem{cor}[thm]{Corollary}
\newtheorem{remark}[thm]{Remark}

\newcommand{\A}{\mathcal{A}}
\newcommand{\CC}{\mathcal{C}}
\newcommand{\QQ}{\mathcal{Q}}
\newcommand{\DD}{\mathcal{D}}
\newcommand{\PP}{\mathcal{P}}
\newcommand{\sd}{D^s}
\newcommand{\bO}{\textbf{O}}
\newcommand{\bP}{\textbf{P}}

\definecolor{myblack}{RGB}{0,0,0}
\definecolor{myred}{RGB}{255,0,0}
\definecolor{myblue}{RGB}{0,0,255}
\definecolor{mygreen}{RGB}{0,120,0}
\definecolor{mypurple}{RGB}{120,0,120}
\definecolor{myorange}{RGB}{255,120,0}
\definecolor{mylightblue}{RGB}{0,255,255}
\newcommand{\red}{\textcolor{myred}{red}}
\newcommand{\blue}{\textcolor{myblue}{blue}}
\newcommand{\green}{\textcolor{mygreen}{green}}
\newcommand{\purple}{\textcolor{mypurple}{purple}}
\newcommand{\orange}{\textcolor{myorange}{orange}}
\newcommand{\lightblue}{\textcolor{mylightblue}{light blue}}

\DeclareMathOperator{\Fun}{Fun}
\DeclareMathOperator{\Hom}{Hom}
\DeclareMathOperator{\sk}{sk}
\DeclareMathOperator{\gSet}{gSet}
\DeclareMathOperator{\id}{id}
\DeclareMathOperator{\Id}{Id}
\DeclareMathOperator{\Map}{Map}
\DeclareMathOperator{\SW}{SW}
\DeclareMathOperator{\Adj}{Adj}
\DeclareMathOperator{\Iso}{Iso}
\DeclareMathOperator{\Cat}{Cat}
\DeclareMathOperator{\h}{h}
\DeclareMathOperator{\I}{I}
\DeclareMathOperator{\II}{II}

\usepackage[nocompress]{cite}
\usepackage[square]{natbib}
\setcitestyle{numbers}

\title{Coherence for adjunctions in a $3$-category via string diagrams}

\author{Manuel Ara\'{u}jo}

\address{Departament of Computer Science and Technology, University of Cambridge, United Kingdom
}
\email{manuel.araujo@cst.cam.ac.uk}

\begin{document}

\setcounter{secnumdepth}{3}
\setcounter{tocdepth}{3}

\maketitle

\begin{abstract}

We construct a $3$-categorical presentation $\Adj_{(3,1)}$ and define a coherent adjunction in a strict $3$-category $\CC$ as a map $\Adj_{(3,1)}\to\CC$. We use string diagrams to show that any adjunction in $\CC$ can be extended to a coherent adjunction in an essentially unique way. The results and their proofs will apply in the context of Gray $3$-categories after the string diagram calculus is shown to hold in that context in an upcoming paper.

\end{abstract}

\section{Introduction}

In this paper, we construct a $3$-categorical presentation $\Adj_{(3,1)}$ containing $1$-cells $l$ and $r$ and we define a \textbf{coherent adjunction} in a strict $3$-category $\CC$ as a functor $\Adj_{(3,1)}\to\CC$. We then prove our Main Theorem, stating that any adjunction in $\CC$ (by which we mean an adjunction in its homotopy $2$-category) can be promoted to a coherent adjunction in an essentially unique way.

In order to state this Theorem precisely, denote by $\theta^{(1)}$ the computad consisting of a single $1$-cell, so that $\Map(\theta^{(1)},\CC)$ is the $3$-groupoid of $1$-morphisms in $\CC$, and let $\Map^L(\theta^{(1)},\CC)$ be its full $3$-subgroupoid whose objects are the left adjoint $1$-morphisms in $\CC$. The map $$E_l:\Map(\Adj_{(3,1)},\CC)\to\Map(\theta^{(1)}, \CC)$$ given by restriction to the $1$-cell $l$ factors through $\Map^L(\theta^{(1)},\CC)$.

\begin{thm}[Main Theorem]

	Given a strict $3$-category $\CC$, the restriction map $$E_l:\Map(\Adj_{(3,1)},\CC)\to\Map^{L}(\theta^{(1)}, \CC)$$ is a weak equivalence of strict $3$-groupoids.

\end{thm}

The basic idea of the proof is that the map $E_l$ is a fibration of $3$-groupoids, by the main result of \cite{adjoints1}. This allows us to make use of a long exact sequence in homotopy groups to reduce the problem to showing that the homotopy groups of the fibre are trivial and the map is surjective on objects. We then prove this is the case by constructing trivialising morphisms for arbitrary elements of these homotopy groups. We do this one cell at a time, by using the lifting properties of fibrations and using the string diagram calculus developed in \cite{Jamie}, \cite{thesis}, \cite{adjoints1} and \cite{sesqui} for explicit constructions.

\begin{remark}

The restriction to strict $3$-categories is a consequence of the fact that we use a string diagram calculus. In an upcoming paper we show that this string diagram calculus holds in Gray $3$-categories and then all results in this paper will hold in that setting, with the same proofs.

\end{remark}

\subsection{Pullbacks and the long exact sequence for a fibration}

In the proof of the main Theorem, we need to make use of the long exact sequence in homotopy groups corresponding to a fibration of $n$-groupoids. We will also need to use pullbacks of maps of $n$-groupoids along a fibration. We therefore state and prove all the necessary results. These should in principle follow from model theoretic arguments in the folk model structure on strict $n$-categories of \cite{folk}. We prefer to give different proofs here for completeness and also because we want proofs that will be applicable in the context of Gray $3$-categories and other types of semistrict $n$-categories.

\subsection{Relation to other work}

We start by defining a $2$-categorical presentation $\Adj_{(2,1)}$ and proving the corresponding coherence result for adjunctions in a $2$-category. This presentation can be deduced directly from the definition of an adjunction as a pair of $1$-morphisms, together with unit and counit $2$-morphisms satisfying two relations, known as the snake relations, or triangle identities. The coherence result in this case is essentially equivalent to the well known result that adjoints are unique up to isomorphism.

The essential difference between $\Adj_{(3,1)}$ and $\Adj_{(2,1)}$ is the appearance of a \textbf{swallowtail relation}, named after the well known singularity. Singularity theory and adjunctions in higher categories are related by the Cobordism Hypothesis (see \cite{lurie}).

The swallowtail relations were fist introduced in the context of adjunctions in \cite{verity_thesis}. There they are part of the definition of a \emph{locally adjoint biadjoint pair} in a \emph{strongly bicategory enriched category} (a kind of semistrict $3$-category). Note that what is called a biadjoint pair in \cite{verity_thesis} is what we call an adjunction in the present paper. In that paper, it is proved that given a biadjoint pair one can modify the triangle isomorphisms in such a way that two additional relations (later called swallowtail relations) are satisfied, yielding a locally adjoint biadjoint pair. A string diagram proof of the analogous result for a strict $3$-category has appeared in \cite{thesis}. There is also a formalized string diagram proof in the proof assistant Globular (see \cite{globular}).

The swallowtail relations also appear in \cite{gurski}, in the more general context of \emph{biadjunctions in tricategories}. There they are depicted in terms of pasting diagrams and the complexity of these diagrams is increased by the presence of various morphisms implementing the weak coherence laws which hold in a tricategory. In that paper, it is proved that any any biequivalence in a tricategory extends to a biadjoint biequivalence, satisfying the swallowtail relations. Again, note that what is called a biequivalence in \cite{gurski} we call here simply an equivalence. What is called a biadjunction in \cite{gurski} we call here an adjunction satisfying the swallowtail relations, i.e. a map $\Adj_{(3,1)}\to \CC$.

In \cite{bruce} the author gives a definition of a \emph{coherent adjoint equivalence} between $2$-categories. This is in particular an adjunction in the $3$-category of $2$-categories, and so the swallowtail relations appear in the definition. The composite $3$-morphisms whose equality is asserted by these relations are depicted as movies of $2$-dimensional string diagrams. This seems to be the first place in the literature where these equations relating the cusp isomorphisms for an adjunction are given the name of swallowtail relations, by analogy with the singularity.

In \cite{piotr} the author proves a coherence result for \emph{duals in monoidal bicategories}. More precisely, they prove that the $2$-groupoid of objects in a monoidal bicategory which admit a dual is equivalent to a $2$-groupoid of \emph{coherent dual pairs}, which can be seen as the $2$-groupoid of maps out of a certain computad, which plays the same role as $\Adj_{(3,1)}$ in the present paper. We can specialize the result in \cite{piotr} to the case of strict monoidal $2$-categories. On the other hand, using the fact that monoidal $2$-categories are just $3$-categories with one object, with duals corresponding to adjoints, we can also specialize the result in the present paper to the context of strict monoidal $2$-categories. The two results on coherence for duals in strict monoidal $2$-categories thus obtained are essentially the same. The main advantage of the methods in the present paper is that the proof is made much simpler by the use of string diagrams, a method which we can currently extend to $4$-categories. Moreover, the result in the present paper does not follow from the one in \cite{piotr}, except in the special case where one considers adjunctions in $3$-category with only one object.

In \cite{riehlverity} the authors construct an $(\infty,2)$-category $\underline{\Adj}$ and prove that the space of functors $\underline{\Adj}\to \Cat_{\infty}$ is equivalent to the space of adjunctions in $\Cat_{\infty}$. Informally, we can think of $\Adj_{(3,1)}$ as a finite presentation for the homotopy $3$-category of $\underline{\Adj}$. This seems to be the first place in the literature where coherence for adjunctions in a higher category $\CC$ is stated in terms of an equivalence of spaces between the space of morphisms in $\CC$ which admit an adjoint and the space of maps into $\CC$ out of a category consisting of a free adjunction. The \emph{strictly undulating squiggles} used there are also a kind of string diagram calculus.

\subsection{Future work}

The use of string diagrams comes with the limitation of applying only to strict $3$-categories. However, we will prove in an upcoming paper that this string diagram calculus is applicable also in Gray $3$-categories and therefore the proofs in \cite{adjoints1} and the present paper also hold more generally. A string diagram calculus for Gray $3$-categories with duals already appears in \cite{surface}.

The methods used in this paper are also extended to dimension $4$ in an upcoming paper, where we construct a $4$-categorical presentation $\Adj_{(4,1)}$ and prove an analogous result for adjunctions in $4$-categories. We will then use this result to give a new proof of the coherence result for fully dualizable objects in a strict symmetric monoidal $3$-category in the author's PhD Thesis \cite{thesis}. The cobordism hypothesis allows us to interpret the corresponding presentation as a finite presentation of the $3$-dimensional fully extended framed bordism category, although this would require the coherence result to be extended to weak symmetric monoidal $3$-categories.

\section{Definitions and basic results}

We now give some necessary definitions and recall the main result from \cite{adjoints1} which we will need in the present paper.

\subsection{Strict $n$-categories}

We think of a strict $n$-category as an algebra over a certain monad $$T_n:\gSet_n\to\gSet_n$$ on the category of $n$-globular sets. This is the monad defined in \cite{OperadsCats}, Chapter 8. Alternatively one can think of a strict $n$-category as a category enriched in strict $(n-1)$-categories with the cartesian product. Given a strict $n$-category $\CC$, we denote by $s,t$ its source and target maps.

\subsection{Equivalences}

In a strict $n$-category, we say that a $k$-morphism $f:x\to y$ is an \textbf{isomorphism} if there exists another $k$-morphism $f:y\to x$ such that $f\circ g=\id_y$ and $g\circ f=\id_x$. We also say that $f$ is \textbf{invertible} and we call $g$ its \textbf{inverse} (one can show that it is unique). However, we are more interested in a weaker version of this, known as \textbf{equivalence}.

\begin{df}

	Let $\CC$ be a strict $n$-category. An $n$-morphism $f:x\to y$ in $\CC$ is an \textbf{equivalence} if it is an isomorphism. When $k<n$, a $k$-morphism $f:x\to y$ in $\CC$ is an \textbf{equivalence} when there is another $k$-morphism $g:y\to x$ and equivalences $f\circ g\to\id_y$ and $g\circ f\to\id_x$ in $\CC$. We say that $x$ is \textbf{equivalent} to $y$, and write $x\simeq y$, if there is an equivalence $x\to y$. When $f:x\to y$ is an equivalence, we also call it \textbf{weakly invertible} and any morphism $g:y\to x$ such that $f\circ g\simeq \id_y$ and $g\circ f\simeq \id_x$ is called a \textbf{weak inverse} to $f$. When $f$ is a $k$-morphism and an equivalence we also call it a \textbf{$k$-equivalence}.

\end{df}

\begin{df}

		An \textbf{$n$-groupoid} is an $n$-category all of whose morphisms are equivalences.

\end{df}

Finally, we use the following notion of weak equivalence for functors, which coincides with the one in the folk model structure of \cite{folk}.

\begin{df}

A functor $F:\CC\to\DD$ between strict $n$-categories is called \textbf{essentially surjective} if for every object $d\in\DD$ there exists an object $c\in\CC$ and an equivalence $F(c)\to d$ in $\DD$. A functor $F:\CC\to\DD$ between strict $n$-categories is called a \textbf{weak equivalence} if it is essentially surjective and for all objects $c_1,c_2\in\CC$ the induced functor $\CC(c_1,c_2)\to\DD(F(c_1),F(c_2))$ is a weak equivalence of $(n-1)$-categories.

\end{df}

\begin{df}

	An $n$-groupoid $G$ is called \textbf{weakly contractible} if the map $G\to *$ is a weak equivalence.

\end{df}

\subsection{Adjunctions}

\begin{df}

An \textbf{adjunction} in a strict $2$-category $\CC$ is a pair of $1$-morphisms $l:X\to Y$ and $r:Y\to X$ together with $2$-morphisms $u:\id_X\to r\circ l$ and $c:l\circ r\to\id_Y$ called the unit and the counit, which satisfy two standard relations, called zigzag, snake or triangle identities.

\end{df}

\begin{df}

Let $\CC$ be a strict $n$-category. We define its \textbf{homotopy $2$-category} to be the strict $2$-category $\h_2(\CC)$ obtained by declaring equivalent $2$-morphisms to be equal.

\end{df}

The following definitions of adjunctions in $n$-categories are adapted from the ones given in \cite{lurie} for the case of $(\infty,n)$-categories.

\begin{df}

An \textbf{adjunction} between $1$-morphisms in a strict $n$-category $\CC$ is a pair of $1$-morphisms $l:X\to Y$ and $r:Y\to X$ together with $2$-morphisms $u:\id_X\to r\circ l$ and $c:l\circ r\to\id_Y$ called the unit and the counit, which determine an adjunction in the homotopy $2$-category $\h_2(\CC)$.

\end{df}

This means that an adjunction in a $3$-category consists of a pair of $1$-morphisms $l:X\to Y$ and $r:Y\to X$ together with unit and counit $2$-morphisms satisfying the usual snake relations or triangle identities up to $3$-isomorphism.

\begin{df}

An adjunction between $k$-morphisms in a strict $n$-category $\CC$ is an adjunction between $1$-morphisms in an appropriate $(n-k+1)$-category of morphisms in $\CC$.

\end{df}

The following Lemma relating equivalences and adjunctions is well known.

\begin{lmm}\label{adjeq}

Let $\CC$ be a strict $n$-category, $f:x\to y$ a $k$-equivalence in $\CC$, $g:y\to x$ a weak inverse and $u:\id_x\to g\circ f$ a $(k+1)$-equivalence. Then there exists a $(k+1)$-equivalence $c:f\circ g\to \id_y$ such that $(f,g,u,c)$ is an adjunction in $\CC$.

\end{lmm}

\begin{proof}

By passing to $\h_2(\Hom(s(x),t(x)))$ we can reduce to the case where $n=2$ and $k=1$. Now we just need to find $c:y\to x$ satisfying the two snake relations. This can be done by using string diagrams, as on the nLab page for adjoint equivalence.
\end{proof}

\subsection{Presentations}

An \textbf{$n$-categorical presentation} is simply a collection of $k$-cells for every $k\leq n+1$, whose sources and targets are composites of lower dimensional cells. We interpret the $(n+1)$-cells as relations. Given an $n$-categorical presentation $\PP$ we denote by $F(\PP)$ the \textbf{$n$-category generated by $\PP$}. Its $k$-morphisms are arbitrary composites of the $k$-cells in $\PP$. Two $n$-morphisms are declared equal when they are related by an $(n+1)$-cell. We sometimes write $\PP\to\CC$ to refer to a functor $F(\PP)\to\CC$.

This can be made precise by using the theory of \textbf{computads}. See \cite{CSPThesis} for a detailed treatment of computads and \cite{adjoints1} for our simplified exposition of how we use them.

We denote by $\theta^{(k)}$ the computad generated by a single $k$-cell, so that functors $\theta^{(k)}\to\CC$ are in canonical bijection with the set of $k$-morphisms in $\CC$.

\subsection{String diagrams}

In \cite{Jamie} a string diagram calculus for $4$-categorical compositions is introduced. The authors introduce the notion of a \textbf{signature}, which consists of sets of generating $k$-cells, for each $k\leq 5$. They then define  a $k$-\textbf{diagram} over a signature, be a $4$-categorical composite of cells. They also introduce \textbf{homotopy generators} which are certain cells encoding coherent versions of the interchange laws that hold in strict $4$-categories. Finally they define a \textbf{signature with homotopy generators} as a signature in which we have specified cells implementing these coherent laws.

In \cite{sesqui} we introduced a monad $T_n^{\sd}$ on globular sets encoding the compositional structure of $n$-dimensional string diagrams. Its algebras are called $n$-sesquicategories and they are $n$-globular sets equipped with strictly associative and unital composition and whiskerig operations, but not satisfying the Godement interchange laws. The notion of a signature then coincides with that of a computad for $T_n^{\sd}$.

In an upcoming paper, we will show how one can define a monad $T_3^{ss}$ by adding to $T_3^{\sd}$ certain operations encoding the homotopy generators. The $T_3^{ss}$-algebras are called semistrict $3$-categories and they are defined precisely so that the string diagram calculus from \cite{Jamie} applies. We also show that they are the same as Gray $3$-categories. We are working on extending this to higher dimensions.

In \cite{adjoints1} we explain how to interpret diagrams over a $4$-signature with homotopy generators as specifying composites in a strict $4$-category, by interpreting the homotopy generators as identity morphisms. In the present paper, we will use these string diagrams in the $n=3$ case, so we include below an informal description of these diagrams.

Given a strict $3$-category $\CC$, we use the string diagram calculus to describe composites of morphisms in any dimension, and to prove identities between composite $3$-morphisms. We read odd dimensional diagrams from left to right and even dimensional diagrams from top to bottom. This means the source of an odd dimensional morphism appears on its left and the source of an even dimensional morphism appears above it.

We denote the composite of two composable $1$-morphisms $f,g$ by the labelled diagram \begin{center}\includegraphics[scale=1.5]{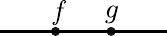} .\end{center} Similarly, we can also denote the composite of $n$ composable $1$-morphisms by a diagram consisting of $n$ labeled dots on a line.

Given $2$-morphisms $\alpha,\beta$ such that $t(\alpha)=s(\beta)$, we can denote their composite by the labeled diagram \begin{center}\includegraphics[scale=1.5]{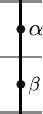} .\end{center} If $f,g$ are $1$-morphisms such that $t(f)=s^2(\alpha)$ and $s(g)=t^2(\alpha)$, we can also denote the whiskering of $\alpha$ with $f$ or $g$ by \begin{center}\begin{tabular}{ccc}\includegraphics[scale=1.5]{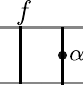} & or & \includegraphics[scale=1.5]{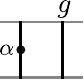} .\end{tabular}\end{center} In general, a diagram such as \begin{center}\includegraphics[scale=1.5]{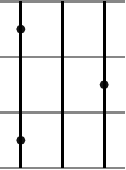} \end{center} labeled by morphisms in $\CC$, subject to compatibility conditions on their sources and targets, determines a composite $2$-morphism in $\CC$. We can also consider $2$-morphisms whose source and target are composites of $1$-morphisms. Given composable $1$-morphisms $i,j$ and $f,g$ we can denote a $2$-morphism $\eta:g\circ f\to j\circ i$ by \begin{center}\includegraphics[scale=1.5]{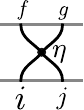} .\end{center} These can also be composed and whiskered with other morphisms, so we can form general $2$-diagrams, which when given a compatible labeling by morphisms in $\CC$ denote composite $2$-morphisms. Here is an example of such a diagram. \begin{center}\includegraphics[scale=1.5]{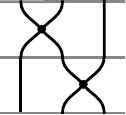} . \end{center}

\newpage 
Now we come to $3$-dimensional diagrams. We denote the composite of two $3$-morphisms by a labeling of \begin{center}\includegraphics[scale=1.5]{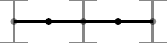} . \end{center} The whiskering of a $3$-morphism with a $2$-morphism corresponds to the diagram \begin{center}\begin{tabular}{ccc}\includegraphics[scale=1.5,align=c]{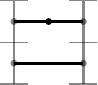} & or & \includegraphics[scale=1.5,align=c]{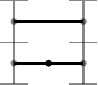} .\end{tabular}\end{center} The whiskering of a $3$-morphism with a $1$-morphism corresponds to the diagram \begin{center}\begin{tabular}{ccc}\includegraphics[scale=1.5,align=c]{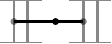} & or & \includegraphics[scale=1.5,align=c]{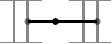} .\end{tabular}\end{center} These basic composition operations can be iterated to form $3$-diagrams such as \begin{center}\includegraphics[scale=1.5]{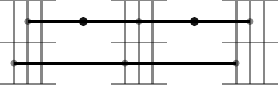} , \end{center} which when given compatible labelings by morphisms in $\CC$ denote composite $3$-morphisms. We can also consider $3$-morphisms whose source and target are arbitrary composites of $1$ and $2$-morphisms in $\CC$, which we denote by labeling a diagram such as \begin{center}\includegraphics[scale=1.5]{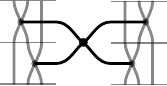} .\end{center} We can then compose these to get general $3$-diagrams, such as \begin{center}\includegraphics[scale=1.5]{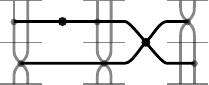} .\end{center}

Notice that the two string diagrams \begin{center}\begin{tabular}{ccc}\includegraphics[scale=1,align=c]{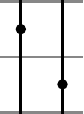} & and & \includegraphics[scale=1,align=c]{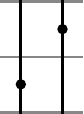}\end{tabular}\end{center} determine the same composition operation on $2$-morphisms in a strict $3$-category, as they both correspond to the pasting diagram $$\xymatrix@1{\bullet\rtwocell & \bullet \rtwocell & \bullet}.$$ So we introduce $3$-dimensional cells \begin{center}\begin{tabular}{ccc}\includegraphics[scale=2,align=c]{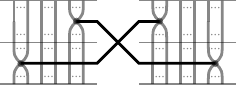} & and & \includegraphics[scale=2,align=c]{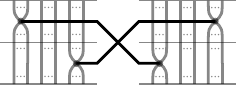} \end{tabular}\end{center} called interchangers (or type $\I_2$ homotopy generators in \cite{Jamie}) which when labeled by compatible morphisms in $\CC$ compose to the appropriate identity $3$-morphism in $\CC$.

\begin{remark}

In a semistrict $3$-category the interchangers become isomorphisms instead of identities.

\end{remark}

Then we need to introduce some equations between $3$-diagrams. First there is interchanger cancellation:

\begin{center}\begin{tabular}{ccc} \includegraphics[scale=2,align=c]{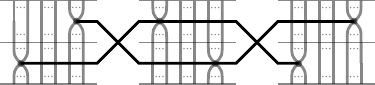} & $=$ &  \includegraphics[scale=2,align=c]{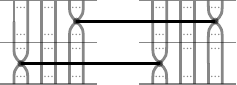} ; \\ \\ \includegraphics[scale=2,align=c]{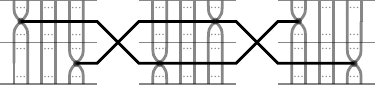} & $=$ &  \includegraphics[scale=2,align=c]{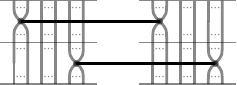} . \\ \\ \end{tabular}\end{center}

Then we have the type $\I_3$ homotopy generator

\begin{center}\begin{tabular}{ccc} \\ \includegraphics[scale=2,align=c]{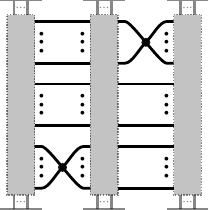} & $=$ &  \includegraphics[scale=2,align=c]{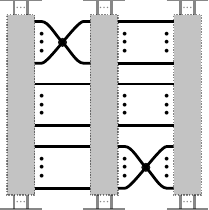} \\ \\  \end{tabular} , \end{center} and finally we have the type $\II_3$ homotopy generators: \begin{center}\begin{tabular}{ccc} \includegraphics[scale=1.8,align=c]{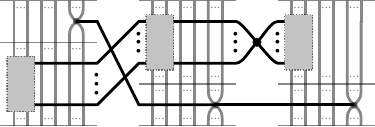} & $=$ &  \includegraphics[scale=1.8,align=c]{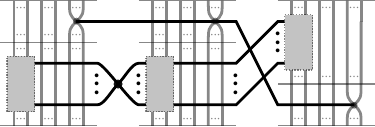} ; \\ \\ \includegraphics[scale=1.8,align=c]{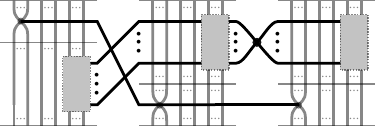} & $=$ &  \includegraphics[scale=1.8,align=c]{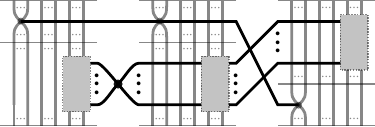} .\\ \\ \end{tabular}\end{center}

\subsection{Functor categories}

Using the left and right internal $\Hom$ from the monoidal biclosed structure on $n$-categories associated to the Crans-Gray tensor product (\cite{crans95}) one can define $n$-categories $\Fun_{lax}(\CC,\DD)$ and $\Fun_{oplax}(\CC,\DD)$ for $n$-categories $\CC$ and $\DD$. One can check that a $k$-morphism in $\Fun_{oplax}(\CC,\DD)$ is a rule that associates to each $\ell$-morphism in $\CC$ a map $\theta^{(k);(\ell)}\to\DD$, satisfying certain relations of compatibility with composition. Here $\theta^{(k);(\ell)}$ is the $(k+\ell)$-computad explictly constructed in \cite{freyd_scheim}. It can also be described as the Crans-Gray tensor product $\theta^{(k)}\otimes\theta^{(\ell)}$. Similarly, a $k$-morphism in $\Fun_{lax}(\CC,\DD)$ is a rule that associates to each $\ell$-morphism in $\CC$ a map $\theta^{(\ell);(k)}\to\DD$.

One can then define the $n$-category $\Fun(\CC,\DD)$ as the subcategory of $\Fun_{oplax}(\CC,\DD)$ consisting of those $k$-morphisms which associate to an $\ell$-morphism in $\CC$ a $(k+\ell)$-equivalence in $\DD$, for $k,\ell\geq 1$. The $k$-morphisms in $\Fun(\CC,\DD)$ are called $k$-transfors. For $k=1,2,3$ they are also called natural transformations, modifications and perturbations, respectively (see the nLab page "transfor" for a discussion of this terminology). Similarly, $\overline{\Fun}(\CC,\DD)$ is the analogous subcategory of $\Fun_{lax}(\CC,\DD)$. Finally $\Map(\CC,\DD)$ and $\overline{\Map}(\CC,\DD)$ are defined as the underlying subgroupoids in $\Fun(\CC,\DD)$ and $\overline{\Fun}(\CC,\DD)$. Given a presentation $\PP$ we write $\Fun(\PP,\DD)$ instead of $\Fun(F(\PP),\DD)$ and similarly for $\Map$.

In \cite{adjoints1} we gave an explicit description of $\Fun(\CC,\DD)$ in terms of string diagrams, when $\CC$ and $\DD$ are $4$-categories. We include here, for convenience, the string diagram description of $k$-transfors between $3$-categories.

\subsubsection{Natural transformations}

Given functors $F,G:\CC\to\DD$, a \textbf{natural transformation}, or $1$-transfor, $\alpha:F\to G$ consists of the following data. We use {\red} and {\blue} to denote the images of objects and morphisms under $F$
and $G$, respectively.
\begin{enumerate}
	\item[0.] For each object $Y\in \CC$ a $1$-morphism $\alpha_Y:F(Y)\to G(Y)$: \begin{center}\begin{tabular}{lcr}$Y=$ \includegraphics[scale=1.5,align=c]{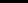} & $\mapsto$ & $\alpha_Y=$  \includegraphics[scale=1.5,align=c]{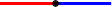}\end{tabular};\end{center}

	\item For each $1$-morphism $g:X\to Y$ in $\CC$ an invertible $2$-morphism $\alpha_g$ in $\DD$: \begin{center}\begin{tabular}{lcr} $g=$ \includegraphics[scale=1.5,align=c]{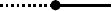} & $\mapsto$ & $\alpha_g=$ \includegraphics[scale=1.5,align=c]{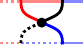} $:$ \includegraphics[scale=1.5,align=c]{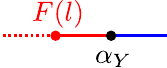} $\to$ \includegraphics[scale=1.5,align=c]{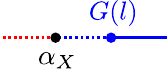}\end{tabular};\end{center}

	\item For each $2$-morphism $\zeta:f\to g$ in $\CC$ an invertible $3$-morphism $\alpha_{\zeta}$ in $\DD$: \begin{center}\begin{tabular}{lcr} $\zeta=$ \includegraphics[scale=1.5,align=c]{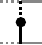} & $\mapsto$ & $\alpha_{\zeta}=$ \includegraphics[scale=1.5,align=c]{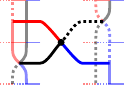} $:$ \includegraphics[scale=1.5,align=c]{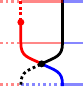} $\to$ \includegraphics[scale=1.5,align=c]{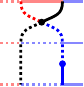}\end{tabular}; \end{center}

	\item For each $3$-morphism $t:\eta\to\zeta$ in $\CC$ a relation $\alpha_{t}$ in $\DD$: \begin{center}\begin{tabular}{lcr}$t=$ \includegraphics[scale=1.5,align=c]{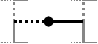} & $\mapsto$ &  \includegraphics[scale=1.5,align=c]{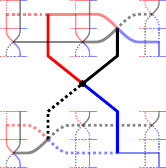} $:$ \includegraphics[scale=1.5,align=c]{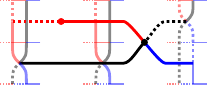} $=$ \includegraphics[scale=1.5,align=c]{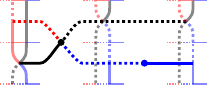}\end{tabular}.\end{center}

\end{enumerate}

This data is subject to relations equating the values of $\alpha$ on composite morphisms with the corresponding composites of values of $\alpha$ given by stacking diagrams.

\subsubsection{Modifications}

Given natural transformations $\alpha,\beta:F\to G$, a \textbf{modification}, or $2$-transfor, $m:\alpha\to\beta$ consists of the following data. We use {\green} for $\alpha$ and {\purple} for $\beta$.

\begin{enumerate}
	\item[0.] For each object $Y\in\CC$ a $2$-morphism $m_Y:\alpha_Y\to \beta_Y$ in $\DD$: \begin{center}\begin{tabular}{lcr} $Y=$ \includegraphics[scale=1.5,align=c]{functorcat/1morph/y.pdf} & $\mapsto$ & $m_Y=$ \includegraphics[scale=1.5,align=c]{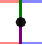} $:$ \includegraphics[scale=1.5,align=c]{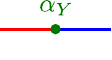} $\to$ \includegraphics[scale=1.5,align=c]{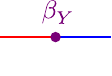}\end{tabular};\end{center}

	\item For each $1$-morphism $g:X\to Y$ in $\CC$ an invertible $3$-morphism $m_g$ in $\DD$: \begin{center}\begin{tabular}{lcr}$g=$ \includegraphics[scale=1.5,align=c]{morphisms/f.pdf} & $\mapsto$ & $m_g=$ \includegraphics[scale=1.5,align=c]{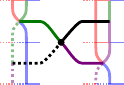} $:$ \includegraphics[scale=1.5,align=c]{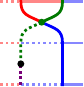} $\to$ \includegraphics[scale=1.5,align=c]{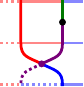};\end{tabular}\end{center}

	\item For each $2$-morphism $\zeta=\includegraphics[scale=1.5,align=c]{morphisms/eta.pdf}:f\to g$ in $\CC$ a relation $m_{\zeta}$ in $\DD$:  \begin{center} \includegraphics[scale=1.5,align=c]{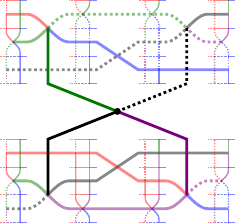}$:$ \includegraphics[scale=1.5,align=c]{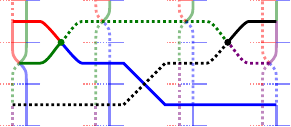}$=$ \includegraphics[scale=1.5,align=c]{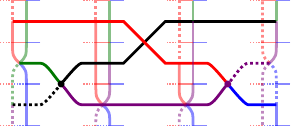} .\end{center}

\end{enumerate}

This data is subject to relations equating the values of $m$ on composite morphisms with the corresponding composites of values of $m$ given by stacking diagrams.

\subsubsection{Perturbations}

Given modifications $l,m:\alpha\to\beta$, a \textbf{perturbation}, or $3$-transfor, $\A:l\to m$ consists of the following data. We use {\orange} for $l$ and {\lightblue} for $m$.

\begin{enumerate}
	\item[0.] For each object $Y\in\CC$ a $3$-morphism $\A_Y:l_Y\to m_Y$ in $\DD$: \begin{center}\begin{tabular}{lcr} $Y=$ \includegraphics[scale=1.5,align=c]{functorcat/1morph/y.pdf} & $\mapsto$ & $\A_Y=$\includegraphics[scale=1.5,align=c]{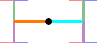} $:$ \includegraphics[scale=1.5,align=c]{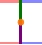} $\to$ \includegraphics[scale=1.5,align=c]{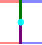}\end{tabular};\end{center}

	\item For each $1$-morphism $g=\includegraphics[scale=1.2,align=c]{morphisms/f.pdf}:X\to Y$ in $\CC$ a relation $\A_g$ in $\DD$: \begin{center}  \includegraphics[scale=1.5,align=c]{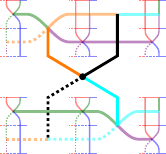} $:$ \includegraphics[scale=1.5,align=c]{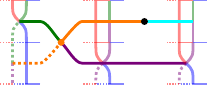} $=$ \includegraphics[scale=1.5,align=c]{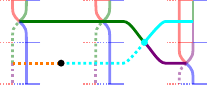} .\end{center}

\end{enumerate}

This data is subject to relations equating the values of $\A$ on composite morphisms with the corresponding composites of values of $\A$ given by stacking diagrams.

\subsection{Fibrations}

\begin{df}

	A map of $n$-groupoids $p:E \to B$ is called a \textbf{fibration} if, given any $k$-morphism $f:x\to y$ in $B$ and a lift $\tilde{x}$ of its source along $p$, there exists a lift $\tilde{f}:\tilde{x}\to\tilde{y}$ of $f$ along $p$.

\end{df}

\begin{remark}

	Given $n$-groupoids $E$ and $B$, it is natural to ask whether a map $f:E\to B$ is a fibration in the sense of this paper if and only if it is is a fibration in the folk model structure on strict $n$-categories defined in \cite{folk}. This is plausible, since the generating trivial cofibrations in that model structure are the inclusions of a free $k$-cell as the source a free fully coherent $(k+1)$-equivalence ($j_k:\bO^k\to \bP^k$ in the notation there). One would therefore have to show that any morphism in an $n$-groupoid can be extended to a fully coherent equivalence. We have not tried to give a proof of this fact.

	Note that in \cite{bg89} the authors construct a model structure on the category of strict $n$-groupoids. However, they define a strict $n$-groupoid as a strict $n$-category where every $k$-morphism has a strict inverse, rather than a weak one. See also \cite{am11}.

\end{remark}

\begin{thm}[from \cite{adjoints1}]\label{fibration}

Let $\CC$ a strict $4$-category, $\PP$ a presentation and $\QQ$ another presentation, obtained by adding a finite number of cells to $\PP$. Then the restriction map $$\Map(\QQ,\CC)\to\Map(\PP,\CC)$$ is a fibration of $4$-groupoids.

\end{thm}

\begin{remark}

In \cite{al20} it is shown that the category of strict $n$-categories equipped with the Crans-Gray tensor product and the folk model structure is a biclosed monoidal model category. This implies that the internal $\Hom$ functors $\Fun_{lax}(-,\DD)$ and $\Fun_{oplax}(-,\DD)$ send cofibrations to fibrations. From \cite{folk} one can deduce that an inclusion of presentations induces a cofibration between the presented $n$-categories. Therefore one can deduce that the restriction map on (op)lax functor categories is a fibration in the folk model structure. Note also that in \cite{al20} it is proved that $\Map(\CC,\DD)$ is the underlying $n$-groupoid of $\Fun_{oplax}(\CC,\DD)$. One might then be able to prove that a folk fibration between lax functor $n$-categories restricts to a fibration (in our sense) between the underlying $n$-groupoids. In this way one might be able to give a different proof of Theorem \ref{fibration} for all $n$. In \cite{adjoints1} we give an explicit string diagram proof of this Theorem in the case $n=4$, which would also apply in any model of semistrict $4$-categories admitting a string diagram calculus.

\end{remark}

\section{Coherence for adjunctions in a $2$-category}

We start by proving coherence for adjunctions in a $2$-category. This is not a new result, as it essentially amounts to uniqueness of adjoints in a $2$-category. We give the proof only to illustrate the general method that will be applied to $3$-categorical and, in a subsequent paper, $4$-categorical adjunctions.

An adjunction in a $2$-category consists of $1$-morphisms $l:X\to Y$ and $r:Y\to X$ together with unit and counit $2$-morphisms satisfying the snake relations. In string diagram notation we can write $l=$ \includegraphics[scale=1,align=c]{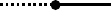} and $r=$ \includegraphics[scale=1,align=c]{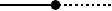}, where we use \includegraphics[scale=2,align=c]{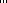} to denote $X$ and \includegraphics[scale=2,align=c]{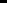} to denote $Y$. If we denote the unit and counit morphisms by \begin{center}\begin{tabular}{lll}\includegraphics[scale=1]{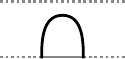} & and & \includegraphics[scale=1]{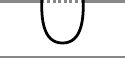}\end{tabular} , \end{center} then the snake relations look like

\begin{center}\begin{tabular}{lllllll}

\includegraphics[scale=1]{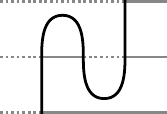} & $=$ & \includegraphics[scale=1]{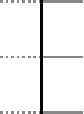} & \text{ and } & \includegraphics[scale=1]{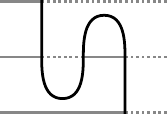} & $=$ & \includegraphics[scale=1]{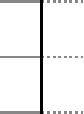} .

\end{tabular}\end{center}

This leads us to make the following definition.

 \begin{df}

 	The presentation $\Adj_{(2,1)}$ consists of

 	\begin{enumerate}
 		\item[0.] $0$-cells $X=$ \includegraphics[scale=2,align=c]{adj/pres/x.pdf} and $Y=$ \includegraphics[scale=2,align=c]{adj/pres/y.pdf} ;
 		\item $1$-cells $l=$ \includegraphics[scale=1,align=c]{adj/pres/l.pdf} $:X\longrightarrow Y$ and $r=$ \includegraphics[scale=1,align=c]{adj/pres/r.pdf} $:Y\longrightarrow X$ ;
 		\item $2$-cells

 			\begin{center}\begin{tabular}{cccccc}

 				$u=$ & \includegraphics[scale=1,align=c]{adj/pres/u.pdf} & $:$ & \includegraphics[scale=1,align=c]{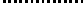} & $\Longrightarrow$ & \includegraphics[scale=1,align=c]{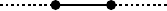} ;

 				\\

 				\\

 				$c=$ & \includegraphics[scale=1,align=c]{adj/pres/c.pdf} & $:$ & \includegraphics[scale=1,align=c]{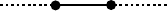} & $\Longrightarrow$ & \includegraphics[scale=1,align=c]{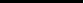} ;

 			\end{tabular}\end{center}

 		\item Relations

 			\begin{center}\begin{tabular}{lllllllll}

 				$C_l=$ & \includegraphics[scale=1,align=c]{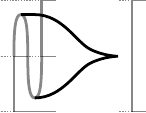} & $:$ & \includegraphics[scale=1,align=c]{adj/pres/snake_l.pdf} & $=$ & \includegraphics[scale=1,align=c]{adj/pres/id_l.pdf} ;

 				\\

 				\\

 				$C_r=$ & \includegraphics[scale=1,align=c]{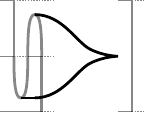} & $:$ & \includegraphics[scale=1,align=c]{adj/pres/snake_r.pdf} & $=$ & \includegraphics[scale=1,align=c]{adj/pres/id_r.pdf} .

 			\end{tabular}\end{center}

 	\end{enumerate}

 \end{df}

\begin{remark}

The $2$-category $F(\Adj_{(2,1)})$ is canonically isomorphic to the $2$-category $\Adj$ defined in \cite{schstr}.

\end{remark}

 Given this definition, we want to prove the following statement.

\begin{p}\label{2d}

	Given a $2$-category $\CC$, the restriction map $$E_l:\Map(\Adj_{(2,1)},\CC)\to\Map^L(\theta^{(1)}, \CC)$$ is a weak equivalence of $2$-groupoids.

\end{p}

%

\begin{df}

Let $f:E\to B$ be a fibration of $2$-groupoids and $b\in B$ an object. The fibre $f^{-1}(b)$ is the $2$-subcategory of $E$ consisting of objects that map to $b$, $1$-morphisms that map to $\id_b$ and $2$-morphisms that map to $\Id_{\id_b}$.

\end{df}

By Theorem \ref{fibration}, $E_l$ is a fibration of $2$-groupoids. In the following sections, we will prove that the fibre of a fibration of $n$-groupoids is an $n$-groupoid and that a fibration which is surjective on objects and has weakly contractible fibres is a weak equivalence. We will also define the homotopy groups $\pi_k(G)$ of an $n$-groupoid $G$ and show that $G$ is weakly contractible if and only if $\pi_k(G)$ is trivial for all $k\leq n$.

\begin{df}

An $n$-groupoid $G$ is \textbf{connected} if for any objects $x,y\in G$ there exists a morphism $x\to y$ in $G$. A connected groupoid $G$ is \textbf{$1$-connected} if given a $1$-morphism $f:x\to x$ in $G$ there exists a $2$-morphism $f\to\id_x$ in $G$.

\end{df}

Once we have defined homotopy groups it will be obvious that $G$ is connected if and only if $\pi_0(G)$ is trivial and $1$-connected when $\pi_0$ and $\pi_1$ are both trivial. In this section, we will show that $E_l$ is surjective on objects and that its fibres are $1$-connected $1$-groupoids, and therefore weakly contractible.

The following is an explicit description of $\Map(\Adj_{(2,1)},\CC)$. An object $F$ in $\Map(\Adj_{(2,1)},\CC)$ is a functor, which consists of a choice of $k$-morphism $F(x)$ for each $k$-cell $x$ in $\Adj_{(2,1)}$, subject to source and target compatibilities. Given functors $F$ and $G$, a $1$-morphism $\alpha:F\to G$ is a weakly invertible natural transformation. Using {\red} and {\blue} to denote the images of generating cells under $F$ and $G$, respectively, $\alpha$ consists of

\begin{enumerate}

 \item[0.] weakly invertible $1$-morphisms \begin{center}$\alpha_X=$ \includegraphics[scale=1,align=c]{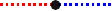} $:F(X)\to G(X)$ and $\alpha_Y=$ \includegraphics[scale=1,align=c]{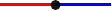} $:F(Y)\to G(Y)$;\end{center}

 \item invertible $2$-morphisms

 \begin{center}\begin{tabular}{llllll}

 $\alpha_l=$ & \includegraphics[scale=1.5,align=c]{functorcat/1morph/a_f.pdf} & $:$ & \includegraphics[scale=1.5,align=c]{functorcat/1morph/a_f_s.pdf} & $\Longrightarrow$ & \includegraphics[scale=1.5,align=c]{functorcat/1morph/a_f_t.pdf} ; \\ \\ $\alpha_r=$ & \includegraphics[scale=1.5,align=c]{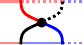} & $:$ & \includegraphics[scale=1.5,align=c]{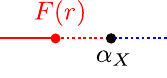} & $\Longrightarrow$ & \includegraphics[scale=1.5,align=c]{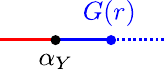};\end{tabular}\end{center}

 \item relations \begin{center}\begin{tabular}{ccc}

 $\alpha_{u} : \includegraphics[scale=1.5,align=c]{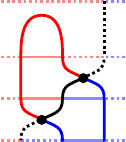} = \includegraphics[scale=1.5,align=c]{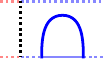}$ & and & $\alpha_{c}: \includegraphics[scale=1.5,align=c]{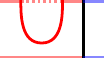} = \includegraphics[scale=1.5,align=c]{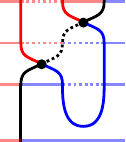}$ . \end{tabular}\end{center}
\end{enumerate}

Finally, given weakly invertible natural transformations $\alpha$, $\beta:F\to G$, a $2$-morphism $\alpha\Rightarrow\beta$ is an invertible modification. Using {\green} and {\purple} to denote the components of $\alpha$ and $\beta$, respectively, $m$ consists of

\begin{enumerate}

	\item[0.] invertible $2$-morphisms \begin{center}\begin{tabular}{ccc}
 $m_X=\includegraphics[scale=1.5,align=c]{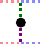}: \includegraphics[scale=1,align=c]{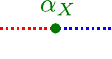} \Longrightarrow \includegraphics[scale=1,align=c]{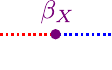}$ & and & $m_Y=\includegraphics[scale=1.5,align=c]{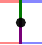} : \includegraphics[scale=1,align=c]{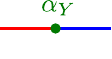} \Longrightarrow \includegraphics[scale=1,align=c]{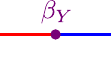}$ ;	\end{tabular}\end{center}

	\item relations \begin{center}\begin{tabular}{ccc}$m_l:\includegraphics[scale=1.5,align=c]{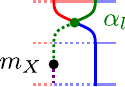} = \includegraphics[scale=1.5,align=c]{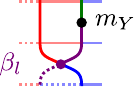}$ & and & $m_r: \includegraphics[scale=1.5,align=c]{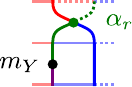}= \includegraphics[scale=1.5,align=c]{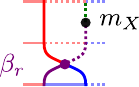}$ .\end{tabular}\end{center}
\end{enumerate}

\begin{lmm}

	Given a $2$-category $\CC$, the restriction map $$E_l:\Map(\Adj_{(2,1)},\CC)\to\Map^L(\theta^{(1)}, \CC)$$ is surjective on objects.

\end{lmm}

\begin{proof}

	Let $F(l):F(X)\to F(Y)$ be a $1$-morphism in $\CC$ which is a left adjoint. We can pick a right adjoint $F(r)$ to $F(l)$, together with unit and counit $2$-morphisms satisfying the snake relations and this data determines a functor $\Adj_{(2,1)}\to\CC$.
\end{proof}

Now we show that the fibres are connected. For this we will need the following Lemmas.

\begin{lmm}\label{ardef}

	Let $\CC$ be a $2$-category. Given $F,G\in\Map(\Adj_{(2,1)},\CC)$ with $F=G$ on $\{X,Y,l\}$ there exists an equivalence $\alpha:F\to G$ in $\Map(\{X,Y,l,r\},\CC)$, which is the identity on $\{X,Y,l\}$.

\end{lmm}

\begin{proof}

	We use red to denote the images of cells under $F$ and blue for their images under $G$. For the images of cells where $F=G$ we use black. So we denote \begin{center}$F(u):=\includegraphics[scale=1,align=c]{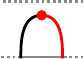}$ , $F(c):=\includegraphics[scale=1,align=c]{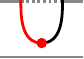}$ , $G(u):=\includegraphics[scale=1,align=c]{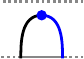}$ and $G(c):=\includegraphics[scale=1,align=c]{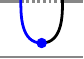}$ .\end{center} We need an isomorphism $\alpha_r:F(r)\to G(r)$, so take \begin{center}$\alpha_r:=$ \includegraphics[scale=1,align=c]{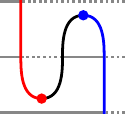} .\end{center}
\end{proof}

\begin{lmm}\label{audef}

	Let $\CC$ be a $2$-category. Given $F,G\in\Map(\Adj_{(2,1)},\CC)$ with $F=G$ on $\{X,Y,l,r\}$ there exists an equivalence $\alpha:F\to G$ in $\Map(\{X,Y,l,r,u\},\CC)$, which is the identity on $\{X,Y,l\}$.

\end{lmm}

\begin{proof}

	We use green to denote the values of $\alpha$. Since $\alpha$ is the identity on $X$, $Y$ and $l$, the relation $\alpha_u$ will be of the form \begin{center}\includegraphics[scale=1.5,align=c]{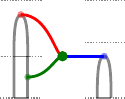} .\end{center} So we define $$\alpha_r=\includegraphics[scale=1,align=c]{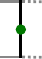}:=\includegraphics[scale=1,align=c]{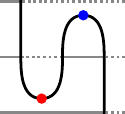}$$ and then the following is a proof of $\alpha_u:$ \begin{center}\includegraphics[scale=1.5,align=c]{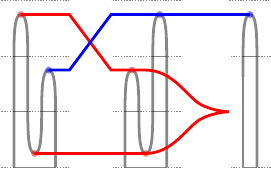} .\end{center}
\end{proof}

\begin{lmm}\label{acdef}

	Let $\CC$ be a $2$-category. Given $F,G\in\Map(\Adj_{(2,1)},\CC)$ with $F=G$ on $\{X,Y,l,r,u\}$ we have $F=G$.

\end{lmm}

\begin{proof}

	We have the following proof that $F(c)=G(c):$ \begin{center}\includegraphics[scale=1]{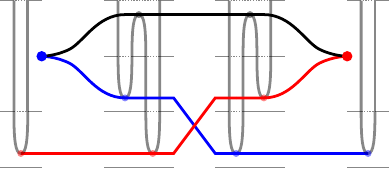} .\end{center}
\end{proof}

\begin{lmm}\label{connected}

	Given a $2$-category $\CC$, the fibres of $E_l:\Map(\Adj_{(2,1)},\CC)\to\Map^L(\theta^{(1)}, \CC)$ are connected.

\end{lmm}

\begin{proof}

	Consider $F,G\in \Map(\Adj_{(2,1)},\CC)$ which agree on $X$ ,$Y$ and $l$. We want to define an equivalence $$\alpha:F\to G$$ in $\Map(\Adj_{(2,1)},\CC)$ which restricts to the identity on $X$, $Y$ and $l$. By Lemma \ref{ardef} there exists an equivalence $\alpha:F\to G$ in $\Map(\{X,Y,l,r\},\CC)$, which is the identity on $\{X,Y,l\}$. Since the restriction map $$\Map(\Adj_{(2,1)},\CC)\to\Map(\{X,Y,l,r\},\CC)$$ is a fibration, one can extend this to an equivalence $\alpha:F\to F_1$ in $\Map(\Adj_{(2,1)},\CC)$, where $F_1$ agrees with $G$ on $X,Y,l,r$ and $\alpha$ is the identity on $X,Y,l$. So now it is enough to find an equivalence $F_1\to G$ which is the identity on $X,Y,l$, where $F_1=G$ on $\{X,Y,l,r\}$. So we can repeat this process, applying the above Lemmas, to get the an equivalence $\alpha:F\to G$ in $\Map(\{X,Y,l,r,u,c\},\CC)$, which is the identity on $\{X,Y,l\}$. Now $\CC$ is a $2$-category and $\{X,Y,l,r,u,c\}$ is the $2$-skeleton of $\Adj_{(2,1)}$, so $1$-morphisms in $\Map(\{X,Y,l,r,u,c\},\CC)$ are the same thing as $1$-morphisms in $\Map(\Adj_{(2,1)},\CC)$.
\end{proof}

Now we show that the fibres are $1$-connected.

\begin{lmm}\label{mrdef}

	Let $\CC$ be a $2$-category. Given a $1$-morphism $\alpha:F\to F$ in $\Map(\Adj_{(2,1)},\CC)$ such that $\alpha$ is the identity on $\{X,Y,l\}$, we have $\alpha=\Id$.

\end{lmm}

\begin{proof}
	Denote $\alpha_r:F(r)\to F(r)$ by \includegraphics[scale=1,align=c]{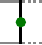} and consider the relation $\alpha_u:$ \begin{center} \includegraphics[scale=1.5,align=c]{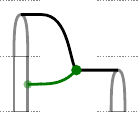} . \end{center} The following is a proof that $\alpha_r=\Id_r:$ \begin{center}\includegraphics[scale=1.5,align=c]{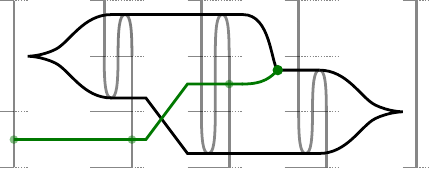} .\end{center}
\end{proof}

\begin{cor}

	Given a $2$-category $\CC$, the fibres of $E_l:\Map(\Adj_{(2,1)},\CC)\to\Map^L(\theta^{(1)}, \CC)$ are $1$-connected.

\end{cor}

\begin{proof}

This follows directly from the previous Lemma.
\end{proof}

\begin{lmm}

	Given a $2$-category $\CC$, the fibres of $E_l:\Map(\Adj_{(2,1)},\CC)\to\Map^L(\theta^{(1)}, \CC)$ are weakly contractible.

\end{lmm}

\begin{proof}

	Since $\CC$ is a $2$-category and $\theta^{(1)}$ contains the $0$-skeleton of $\Adj_{(2,1)}$, the fibres of this map are $1$-groupoids. Therefore, since the fibres are $1$-connected, they are weakly contractible.
\end{proof}

\section{The homotopy pullback}

Given a diagram of $n$-groupoids $$\xymatrix{ & X\ar[d]^{F} \\ Y\ar[r]_{G} & Z,}$$ its pullback $X\times_{Z} Y$ in $\Cat_n$ is just the pullback in $\gSet_n$, equipped with an obvious $T_n$-algebra structure. We want to show that $X\times_{Z} Y$ is actually an $n$-groupoid, provided either $F$ or $G$ is a fibration. One could give a direct proof of this, using the lifting properties of the fibration and the fact that equivalences can be promoted to adjoint equivalences to construct weak inverses. However, we prefer to give another proof, using computads and Theorem \ref{fibration}, which is more in tune with the general idea of this paper. Our strategy is to define a homotopy pullback $X\times_{Z}^h Y$, which we can show is always an $n$-groupoid, and then to show that when $F$ or $G$ is a fibration the natural map $$X\times_{Z} Y\to X\times_{Z}^h Y$$ is a weak equivalence of $n$-categories.

\begin{df}

We define the \textbf{homotopy pullback} of a diagram $$\xymatrix{ & X\ar[d]^{F} \\ Y\ar[r]_{G} & Z}$$ of $n$-groupoids to be the pullback $$\xymatrix{X\times_{Z}^h Y\ar[r]\ar[d] & X\times Y\ar[d]^{F\times G} \\ \overline{\Map}(\theta^{(1)},Z)\ar[r]_-{s\times t} & Z\times Z }$$ in $\Cat_n$.

\end{df}

\begin{remark}

Given an $n$-groupoid $Z$, a $k$-morphism in $\overline{\Map}(\theta^{(\ell)},Z)$ is the same thing as an $\ell$-morphism in $\Map(\theta^{(k)},Z)$, since they both correspond to maps $\theta^{(\ell);(k)}\to Z$.

\end{remark}

So the set of $k$-morphisms in the homotopy pullback is the pullback of sets

$$\xymatrix{\Hom(\theta^{(k)},X\times_{Z}^h Y)\ar[r]\ar[d] & \Hom(\theta^{(k)},X)\times\Hom(\theta^{(k)},Y)\ar[d] \\ \Hom(\theta^{(k)},\overline{\Map}(\theta^{(1)},Z))\ar[r] & \Hom(\theta^{(k)},Z\times Z),}$$ which is equal to the pullback of sets

$$\xymatrix{\Hom(\theta^{(k)},X\times_{Z}^h Y)\ar[r]\ar[d] & \Hom(\theta^{(k)},X)\times\Hom(\theta^{(k)},Y)\ar[d] \\ \Hom(\theta^{(1)},\Map(\theta^{(k)},Z))\ar[r] & \Hom(\partial \theta^{(1)},\Map(\theta^{(k)},Z)).}$$ This means that a $k$-morphism in $X\times_{Z}^h Y$ consists of a diagram $$\xymatrix{\theta^{(k)}\ar[r]\ar[d] & X\ar@2[dl]\ar[d] \\ Y\ar[r] & Z.}$$

Now we show that the homotopy pullback of $n$-groupoids is an $n$-groupoid.

\begin{df}

Let $\Iso_{(k,k)}$ be the computad with two parallel $m$-cells $x_m$ and $y_m$ for each $m<k$, two $k$-cells $l:x_{k-1}\to y_{k-1}$ and $r:y_{k-1}\to x_{k-1}$ and two $(k+1)$-cells $u:\id_{x_{k-1}}\to r\circ l$ and $c:l\circ r\to\id_{y_{k-1}}$.

\end{df}

\begin{df}

Let $\Adj_{(k+1,k)}$ be the computad obtained from $\Iso_{(k,k)}$ by adding two $(k+2)$-cells, corresponding to the two snake relations.
\end{df}

\begin{lmm}

Let $\CC$ be an $n$-category and suppose all $(k+1)$-morphisms are weakly invertible in $\CC$. Then a $k$-morphism in $\CC$ is weakly invertible if and only if the corresponding map $\theta^{(k)}\to\CC$ extends to $\Iso_{(k,k)}\to\CC$.

\end{lmm}

\begin{proof}

Given a $k$-morphism $l:x\to y$ in $\CC$, an extension of the corresponding map $\theta^{(k)}\to\CC$ to $\Iso_{(k,k)}\to\CC$ consists of a choice of $k$-morphism $r:y\to x$ in $\CC$ and $(k+1)$-morphisms $\id_{x}\to r\circ l$ and $l\circ r\to\id_{y}$. Since all $(k+1)$-morphisms are weakly invertible in $\CC$, such choices exist if and only if $l$ is weakly invertible.
\end{proof}

\begin{lmm}

Let $X$ be an $n$-groupoid and consider $F,G\in\Map(\Adj_{(k+1,k)},X)$ where $F(l)=G(l)$. Then there exists an equivalence $\alpha:F\to G$ in $\Map(\Iso_{(k,k)},X)$, restricting to the identity on $l$.

\end{lmm}

\begin{proof}

 We need to construct $\alpha_r$, $\alpha_{u}$ and $\alpha_{c}$ in $X$. By passing to $\h_2(\Hom_X(F(x_{k-2}),F(y_{k-2})))$ we can reduce to the case where $n=2$ and $k=1$ and apply Lemma \ref{connected}. From this we get $\alpha_r$ as a $2$-morphism and $\alpha_{u}$, $\alpha_{c}$ as identities between $2$-morphisms in $\h_2(\Hom_X(F(x_{k-2}),F(y_{k-2})))$. These correspond to the required $(k+1)$-morphism $\alpha_r$ and the $(k+2)$-morphisms $\alpha_{u}$ and $\alpha_{c}$ in $X$.
\end{proof}

\begin{lmm}

Let $X$ be an $n$-groupoid, $F,G\in\Map(\Adj_{(k+1,k)},X)$ functors and $\alpha:F\to G$ in $\Map(\theta^{(k)},X)$ an equivalence. Then $\alpha$ extends to an equivalence $F\to G$ in $\Map(\Iso_{(k,k)},X)$.

\end{lmm}

\begin{proof}

Passing to a $\Hom_X(F(x_{k-2}),F(y_{k-2}))$ we can reduce to the case where $k=1$. Now $\Adj_{(2,1)}$ only has cells of dimension $\leq 3$, so extending $\alpha$ only involves constructing composites of dimension $\leq 4$ in $\Hom_X(F(x_{k-2}),F(y_{k-2}))$.

Therefore we can apply Theorem \ref{fibration} to $\h_3(\Hom_X(F(x_{k-2}),F(y_{k-2})))$ and lift $\alpha$ starting at $F$ to get an equivalence $F\to G_1$ in $\Map(\Adj_{(k+1,k)},X)$, where $G_1(l)=G(l)$. Now we just need to find an equivalence $G_1\to G$ in $\Map(\Iso_{(k,k)},X)$ restricting to the identity on $l$, which we can do by the previous Lemma.
\end{proof}

\begin{p}

The homotopy pullback of a diagram of $n$-groupoids is an $n$-groupoid.

\end{p}

\begin{proof}

Consider a diagram $$\xymatrix{ & X\ar[d]^{F} \\ Y\ar[r]_{G} & Z}$$ of $n$-groupoids. Let $k\leq n$ and suppose all $(k+1)$-morphisms are weakly invertible in the homotopy pullback. Note that this condition vacuously holds for $k=n$. Let $$\xymatrix{\theta^{(k)}\ar[r]^{x}\ar[d]_{y} & X\ar@2[dl]_{\alpha}\ar[d]^{F} \\ Y\ar[r]_{G} & Z}$$ be a $k$-morphism in the homotopy pullback. We want to show that it is weakly invertible, so we need to find an extension

 $$\xymatrix{\theta^{(k)}\ar[dr]\ar@/^/[drr]^{x}\ar@/_/[ddr]_{y} & & \\ &\Iso_{(k,k)}\ar@{.>}[r]\ar@{.>}[d] & X\ar@{:>}[dl]\ar[d]^{F} \\ & Y\ar[r]_{G} & Z.}$$ Since $X$ and $Y$ are $n$-groupoids, Lemma \ref{adjeq} tells us we can find extensions of $x$ and $y$ to $\Adj_{(k+1,k)}$, to get a diagram of the form $$\xymatrix{\theta^{(k)}\ar[dr]\ar@/^/[drr]^{x}\ar@/_/[ddr]_{y} & & \\ &\Adj_{(k+1,k)}\ar[r]^{x}\ar[d]_{y} & X \\ & Y. & }$$

The previous Lemma now allows us to extend $\alpha:F\circ x\to G\circ y$ from $\Map(\theta^{(k)},Z)$ to $\Map(\Iso_{(k,k)},Z)$ as desired.
\end{proof}

\begin{p}

Given a diagram of $n$-groupoids $$\xymatrix{ & X\ar[d]^{F} \\ Y\ar[r]_{G} & Z,}$$ where $F$ is a fibration, the canonical map $$X\times_{Z} Y\to X\times_{Z}^h Y$$ is a weak equivalence of $n$-categories.

\end{p}

\begin{proof}

Take a $k$-morphism $$\xymatrix{\theta^{(k)}\ar[r]^{x}\ar[d]_{y} & X\ar@2[dl]_{\alpha}\ar[d]^{F} \\ Y\ar[r]_{G} & Z.}$$ in the homotopy pullback, whose source and target are in the image of the map $X\times_{Z} Y\to X\times_{Z}^h Y$. This means that $\alpha_{\partial \theta^{(k)}}$ is the identity natural transformation and so $\alpha$ is simply an equivalence $F(x)\to G(y)$ in $Z$. Since $F$ is a fibration, we can find a lift $\phi$ in

$$\xymatrix{\theta^{(k)}\ar[r]^{x}\ar[d]_{s} & X \ar[d]^{F} \\ \theta^{(k+1)}\ar@{.>}[ru]^{\phi}\ar[r]_{\alpha} & Z.}$$ Now $\phi$ is a $(k+1)$-morphism $x\to\bar{x}$ such that $F(\phi)=\alpha:F(x)\to G(y)$, so $F(\bar{x})=G(y)$ and we have a $k$-morphism $$\xymatrix{\theta^{(k)}\ar[r]^{\bar{x}}\ar[d]_{y} & X\ar@{=}[dl]\ar[d]^{F} \\ Y\ar[r]_{G} & Z}$$ in the pullback $X\times_{Z} Y$. In order to show that the canonical map is essentially surjective on $k$-morphisms, we need to show that this is equivalent to the original $k$-morphism in the homotopy pullback. Since $X\times_{Z}^h Y$ is an $n$-groupoid, it's enough to show that there is a $(k+1)$-morphism between them. So we need to find a natural transformation $$\xymatrix{\theta^{(k+1)}\ar[r]^{\phi}\ar[d]_{\id_y} & X\ar@{:>}[dl]\ar[d]^{F} \\ Y\ar[r]_{G} & Z}$$ whose restrictions along $\xymatrix{\theta^{(k)}\ar[r]^{s} & \theta^{(k+1)}}$ and $\xymatrix{\theta^{(k)}\ar[r]^{t} & \theta^{(k+1)}}$ are

 $$\xymatrix{\theta^{(k)}\ar[r]^{x}\ar[d]_{y} & X\ar@2[dl]_{\alpha}\ar[d]^{F} & & \theta^{(k)}\ar[r]^{\bar{x}}\ar[d]_{y} & X\ar@{=}[dl]\ar[d]^{F} & \\ Y\ar[r]_{G} & Z & \text{and} & Y\ar[r]_{G} & Z & \text{respectively.}}$$ We can extend to the $(k+1)$-cell by $$\xymatrix{F(x)\ar[r]^{F(\phi)}\ar[d]_{\alpha} & F(\bar{x})\ar@{=}[dl]\ar[d]^{\id} \\ G(y)\ar[r]_{\id} & G(y).}$$

\end{proof}

\begin{lmm}

Consider a weak equivalence of $n$-categories $\CC\to\DD$ and suppose that $\DD$ is an $n$-groupoid. Then $\CC$ is an $n$-groupoid.

\end{lmm}

\begin{proof}

Consider a $k$-morphim $f:x\to y$ in $\CC$ and let $F(f)^{-1}:F(y)\to F(x)$ be a weak inverse for $F(f)$. Since $F$ is a weak equivalence, there exists an $\overline{f}:y\to x$ in $\CC$ with $F(\overline{f})\simeq F(f)^{-1}$. Then $F(\overline{f}\circ f)\simeq F(f)^{-1}\circ F(f) \simeq \id_{F(x)}$. Since $F$ is a weak equivalence, this implies $\overline{f}\circ f\simeq \id_{x}$. Similarly, we have $f\circ\overline{f}\simeq \id_{y}$, so $\overline{f}$ is weak inverse to $f$.
\end{proof}

\begin{p}\label{pbngrpd}

Given a diagram of $n$-groupoids $$\xymatrix{ & X\ar[d]^{F} \\ Y\ar[r]_{G} & Z,}$$ where $F$ is a fibration, the pullback $X\times_{Z}Y$ is an $n$-groupoid.

\end{p}

\begin{proof}

This follows from the fact that $X\times_{Z}^h Y$ is an $n$-groupoid and the canonical map $X\times_{Z}Y\to X\times_{Z}^h Y$ is a weak equivalence.
\end{proof}

\section{The long exact sequence for a fibration}

Now we show that a fibration of $n$-groupoids is a weak equivalence if and only if its fibres are weakly contractible, by using an analog of the long exact sequence in homotopy groups for a fibration of spaces.

\begin{lmm}\label{homfib}

Let $p:E\to B$ be a fibration of $n$-groupoids. Then, for any $0\leq k\leq n-1$ and any $k$-morphisms $x,y$ in $E$, the induced map $$E(x,y)\to B(p(x),p(y))$$ is a fibration of $(n-k-1)$-groupoids.

\end{lmm}

\begin{proof}

The lifting condition for $\ell$-morphisms in $B(p(x),p(y))$ follows easily from the lifting condition for $(\ell+k+1)$-morphisms in $B$.
\end{proof}

\begin{df}

Given a map of $n$-groupoids $f:A\to B$ and an object $b\in B$, we define the \textbf{fibre} $f^{-1}(b)$ of $f$ over $b$ to be the pullback

$$\xymatrix{f^{-1}(b)\ar[r]\ar[d] & A\ar[d]^{f} \\ \theta^{(0)} \ar[r]_{b} & B}$$ in $\Cat_n$.

\end{df}

So a $k$-morphism in $f^{-1}(b)$ is a $k$-morphism $a\in A_k$ such that $f(a)=\Id^{(k)}_b$ (defined inductively by saying that $\Id^{(k)}_b$ is the identity morphism on $\Id^{(k-1)}_b$).

\begin{lmm}

If $p:E\to B$ is a fibration between $n$-groupoids and $b\in B$ is an object, then $p^{-1}(b)$ is an $n$-groupoid.

\end{lmm}

\begin{proof}

This follows from Proposition \ref{pbngrpd}.
\end{proof}

\begin{df}

Let $G$ be an $n$-groupoid and $x\in G$ an object. Define $\id^{(0)}_x:= x$ and $\id^{(k)}_x=\id_{\id^{(k-1)}_x}$. Denote by $\Omega_x G$ the $(n-1)$-groupoid $\Hom(x,x)$. This comes equipped with a strictly associative monoidal structure, given by composition in $G$. Denote $\Omega^k_xG:=\Omega_{\id^{(k-1)}_x}\Omega^{k-1}_xG$. Finally, denote by $G_0$ the set of objects of $G$.

\end{df}

\begin{df}

Let $G$ be an $n$-groupoid. We define $\pi_0(G):=G_0/\sim$, where the equivalence relation $\sim$ is equivalence in $G$. Now let $x\in G$ be an object. Define $\pi_0(G,x)$ to be the pointed set $(\pi_0(G),[x])$, were $[x]$ denotes the equivalence class of $x$ in $\pi_0(G)$. Finally, for $1\leq k\leq n$, define $\pi_k(G,x):=\pi_{k-1}(\Omega_xG,\id_x)$ with monoid structure induced by composition.

\end{df}

Note that, for $k\geq 1$, the monoids $\pi_k(G,x)$ are actually groups and for $k\geq 2$ they are abelian, by an Eckmann-Hilton argument with pasting diagrams. Moreover, given a map of $n$-groupoids $f:A\to B$ and an object $a\in A$ one can also define $\pi_k(f,a):\pi_k(A,a)\to\pi_k(B,f(a))$, making $\pi_k$ into a functor on pointed $n$-groupoids.

\begin{lmm}\label{htpygrps}

Let $f:A\to B$ be a map of $n$-groupoids. Then $f$ is a weak equivalence if and only if the maps $\pi_k(f,a):\pi_k(A,a)\to\pi_k(B,f(a))$ are isomorphisms, for all $a\in A_0$ and for all $k\geq 0$.

\end{lmm}

\begin{proof}

The proof is by induction on $n$. Suppose $f$ is a weak equivalence. Then it is essentially surjective, so it is surjective on $\pi_0$. Moreover, for any objects $x,y\in A$, the map $A(x,y)\to B(f(x),f(y))$ is a weak equivalence of $(n-1)$-groupoids, so by the induction hypothesis it induces isomorphisms on all homotopy groups. In particular, it induces isomorphisms $\pi_k(A(x,x),\id_x)\to \pi_k(B(f(x),f(x)),\id_{f(x)})$, for $0\leq k \leq n-1$. Now
\[\pi_k(A(x,x),\id_x)=\pi_{k+1}(A,x)\quad \text{and}\quad \pi_k(B(f(x),f(x)),\id_{f(x)})=\pi_{k+1}(B,f(x)), \]
so we conclude that $\pi_l(f,x):\pi_l(A,x)\to\pi_l(B,f(x))$ is an isomorphism, for $1\leq l \leq n$. So all that is left to do is to show that $\pi_0(f)$ is injective. So suppose $[f(x)]=[f(y)]$ in $\pi_0(B)$ and pick an equivalence $\beta:f(x)\to f(y)$ in $B$. Since $A(x,y)\to B(f(x),f(y))$ is essentially surjective, there exists $\alpha:x\to y$ and an equivalence $f(\alpha)\Rightarrow \beta$. In particular, we have $[x]=[y]$ in $\pi_0(A)$.

Now suppose the maps $\pi_k(f,a):\pi_k(A,a)\to\pi_k(B,f(a))$ are isomorphisms, for all objects $a\in A$ and for all $k\geq 0$. Then $f$ is essentially surjective, being an isomorphism on $\pi_0$. Now let $x,y\in A$ and consider the map $A(x,y)\to B(f(x),f(y))$. We need to show that this map is a weak equivalence of $(n-1)$-groupoids, and by the induction hypothesis it is enough to show that the maps $$\pi_k(A(x,y),\alpha)\to\pi_k(B(f(x),f(y)),f(\alpha))$$ are isomorphisms, for $0\leq k \leq n-1$ and $\alpha:x\to y$ in $A$. Consider the maps $A(x,x)\to A(x,y)$ and $B(f(x),f(x))\to B(f(x),f(y))$ defined by wiskering with $\alpha$ and $f(\alpha)$, respectively. These induce isomorphisms $$\pi_k(A(x,x),\id_x)\to\pi_k(A(x,y),\alpha)$$ and $$\pi_k(B(f(x),f(x)),\id_{f(x)})\to\pi_k(B(f(x),f(y)),f(\alpha))$$ for $0\leq k \leq n-1$, which fit in a commuting diagram

$$\xymatrix{\pi_k(A(x,y),\alpha) \ar[r] & \pi_k(B(f(x),f(y)),f(\alpha)) \\
\pi_k(A(x,x),\id_x)\ar[u]^*[@]-{\simeq} \ar[r]^-{\simeq} & \pi_k(B(f(x),f(x)),\id_{f(x)})\ar[u]^*[@]-{\simeq} \\
\pi_{k+1}(A,x)\ar[r]^{\simeq}\ar@{=}[u] & \pi_{k+1}(B,f(x)).\ar@{=}[u]
}$$This shows that $\pi_k(A(x,y),\alpha)\to\pi_k(B(f(x),f(y)),f(\alpha))$ is an isomorphism.
\end{proof}

Now we construct the long exact sequence in homotopy groups associated with a fibration of $n$-groupoids.

\begin{lmm}

Let $p:E\to B$ be a fibration of $n$-groupoids, $x\in E$ an object, $b=p(x)\in E$ and let $F_b$ be the fibre of $p$ over $b$. Then the induced map $\Omega_xE\to\Omega_bB$ is a fibration of $(n-1)$-groupoids, whose fibre over $\id_b$ is $\Omega_x F_b$.

\end{lmm}

\begin{proof}

The fact that $\Omega_x E\to\Omega_b B$ is a fibration follows from Lemma \ref{homfib}. The fact that the fibre is $\Omega_x F_b$ follows from unraveling the definitions.
\end{proof}

\begin{df}

Let $p:E\to B$ be a fibration of $n$-groupoids, $x\in E$ an object, $b=p(x)\in E$ and let $F_b$ be the fibre of $p$ over $b$. We define a pointed map $\partial:\pi_0(\Omega_b B,\id_b)\to\pi_0 (F_b,x)$ sending the class $[\alpha]$ of a $1$-morphism $\alpha:b\to b$ in $B$ to the class $[y]$ where $y$ is the endpoint of a lift of $\alpha$ along $p$, starting at $x$.

\end{df}

\begin{lmm}

The above procedure gives a well defined pointed map $$\partial:\pi_0(\Omega_b B,\id_b)\to\pi_0 (F_b,x).$$

\end{lmm}

\begin{proof}

Suppose we have $1$-morphisms $f,g:b\to b$ in $B$, and a $2$-morphism $\alpha:f\to g$, so that $[f]=[g]$ in $\pi_0(\Omega_bB,\id_b)$. Let $\bar{f}:x\to y$ and $\bar{g}:x\to z$ be lifts of $f$ and $g$ along $p$. We want to show that $[y]=[z]$ in $F_b$, so we need to find a morphism $y\to z$ in $E$ that maps to $\id_b$. Let $\bar{f}^{-1}$ be a weak inverse for $\bar{f}$ and let $f^{-1}:=p(\bar{f}^{-1})$. Consider the map $\bar{g}\circ\bar{f}^{-1}:y\to z$. This maps down to $g\circ f^{-1}$. Composing $\alpha^{-1}$ with a $2$-morphism $f\circ f^{-1}\to\id_b$, we get a $2$-morphism $g\circ f^{-1}\to\id_b$. Now we can lift this $2$-morphism starting at $\bar{g}\circ \bar{f}^{-1}$ and the target of the resulting $2$-morphism is a $1$-morphism $y\to z$ which maps to $\id_b$ under $p$.
\end{proof}

\begin{lmm}

Let $p:E\to B$ be a fibration of $n$-groupoids and $F_b$ its fibre over $b\in B$. Then the map $\xymatrix{\pi_k(B,b)\ar[r]^{\partial} & \pi_{k-1}(F_b,x)}$ is a group homomorphism for $k\geq 2$.

\end{lmm}

\begin{proof}

It is enough to show that $\xymatrix{\pi_0(\Omega^2_b B,\Id_{\id_b})\ar[r]^{\partial} & \pi_{0}(\Omega_xF_b,\id_x)}$ is a group homomorphism. This follows from the fact that $$\xymatrix@1{b\ruppertwocell^{\id_b}\rlowertwocell_{\id_b}\ar[r] & b} = \xymatrix@1{b\rtwocell^{\id_b}_{\id_b} & b \rtwocell^{\id_b}_{\id_b} & b}.$$
\end{proof}

\begin{p}

Let $p:E\to B$ be a fibration, $x\in E$ an object, set $b:=p(x)$ and let $F_b$ be the fibre of $f$ over $b$. Then the following is an exact sequence of groups and pointed sets

$$\xymatrix{\cdots\ar[r] & \pi_1(E,x)\ar[r] & \pi_1(B,b)\ar[r]^{\partial} & \pi_0(F_b,x)\ar[r] & \pi_0(E,x)\ar[r] & \pi_0(B,b)}.$$

\end{p}

\begin{proof}

It is enough to show that $$\xymatrix{\pi_0(\Omega_xE,\id_x)\ar[r] & \pi_0(\Omega_bB,\id_b)\ar[r]^{\partial} & \pi_0(F_b,x)\ar[r] & \pi_0(E,x)\ar[r] & \pi_0(B,b)}$$ is exact. This in turn is easy to check directly.
\end{proof}

\begin{cor}

A fibration of $n$-groupoids $f:A\to B$ is a weak equivalence if and only if for every object $b\in B$ the fibre $f^{-1}(b)$ is weakly contractible.

\end{cor}

\section{Coherence for adjunctions in a $3$-category}

We now define the presentation $\Adj_{(3,1)}$ consisting of coherence data for an adjunction between $1$-morphisms in a $3$-category and prove our Main Theorem, which we restate here for convenience.

\begin{thm}\label{main}

	Given a strict $3$-category $\CC$, the restriction map $$E_l:\Map(\Adj_{(3,1)},\CC)\to\Map^{L}(\theta^{(1)}, \CC)$$ is a weak equivalence of strict $3$-groupoids.

\end{thm}

 \begin{df}

The presentation $\Adj_{(3,1)}$ consists of

\begin{enumerate}

 \item[0.] objects $X=$ \includegraphics[scale=2,align=c]{adj/pres/x.pdf} and $Y=$ \includegraphics[scale=2,align=c]{adj/pres/y.pdf};\\

 \item $1$-cells $l$= \includegraphics[scale=1,align=c]{adj/pres/l.pdf} $:\xymatrix@1{X\ar@1@<1ex>[r] & \ar@1@<1ex>[l]Y}:$ \includegraphics[scale=1,align=c]{adj/pres/r.pdf} $=r$;\\

 \item $2$-cells

 \begin{center}\begin{tabular}{llllll}

 $u=$ & \includegraphics[scale=1,align=c]{adj/pres/u.pdf} & $:$ & \includegraphics[scale=1,align=c]{adj/pres/u_s.pdf} & $\xymatrix@1{\ar@2[r] & }$ & \includegraphics[scale=1,align=c]{adj/pres/u_t.pdf} ;

 \\

 \\

 $c=$ & \includegraphics[scale=1,align=c]{adj/pres/c.pdf} & $:$ & \includegraphics[scale=1,align=c]{adj/pres/c_s.pdf} & $\xymatrix@1{\ar@2[r] & }$ & \includegraphics[scale=1,align=c]{adj/pres/c_t.pdf} ;

 \end{tabular}\end{center}

 \item $3$-cells

 \begin{center}\begin{tabular}{lllllllll}

 $C_l=$ & \includegraphics[scale=1,align=c]{adj/pres/cusp_l.pdf} & $:$ & \includegraphics[scale=1,align=c]{adj/pres/snake_l.pdf} & $\xymatrix@1{\ar@3@<1ex>[r] & \ar@3@<1ex>[l]}$ & \includegraphics[scale=1,align=c]{adj/pres/id_l.pdf} & $:$ & \includegraphics[scale=1,align=c]{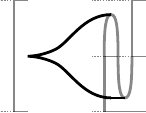} & $=C_l^{-1}$ ;

 \\

 \\

 $C_r=$ & \includegraphics[scale=1,align=c]{adj/pres/cusp_r.pdf} & $:$ & \includegraphics[scale=1,align=c]{adj/pres/snake_r.pdf} & $\xymatrix@1{\ar@3@<1ex>[r] & \ar@3@<1ex>[l]}$ & \includegraphics[scale=1,align=c]{adj/pres/id_r.pdf} & $:$ & \includegraphics[scale=1,align=c]{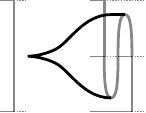} & $=C_r^{-1}$ ;

 \end{tabular}\end{center}

\newpage
 \item relations

 \begin{center}\begin{tabular}{lllll}

 \includegraphics[scale=1,align=c]{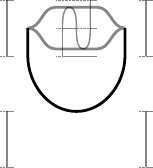} & $:$ & \includegraphics[scale=1,align=c]{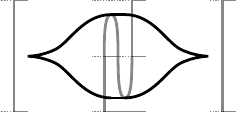} & $=$ & $\Id^{(2)}_l$ ;

 \\

 \\

 \includegraphics[scale=1,align=c]{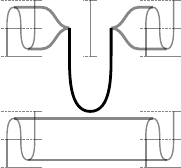} & $:$ & \includegraphics[scale=1,align=c]{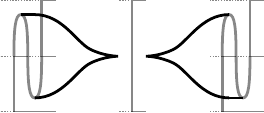} & $=$ & \includegraphics[scale=1,align=c]{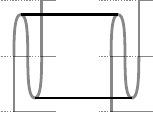} ;

 \\

 \\

 \includegraphics[scale=1,align=c]{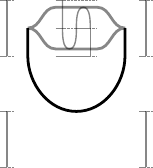} & $:$ & \includegraphics[scale=1,align=c]{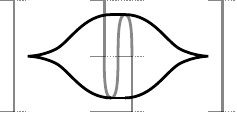} & $=$ & $\Id^{(2)}_r$ ;

 \\

 \\

 \includegraphics[scale=1,align=c]{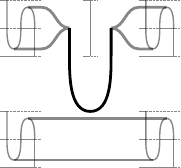} & $:$ & \includegraphics[scale=1,align=c]{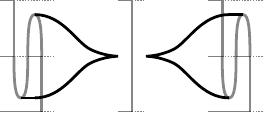} & $=$ & \includegraphics[scale=1,align=c]{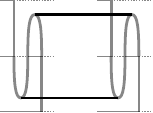} ;

 \\

 \\

 \includegraphics[scale=1,align=c]{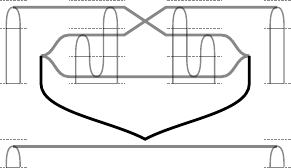} & $:$ & \includegraphics[scale=1,align=c]{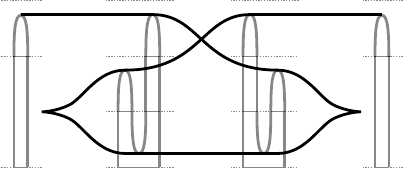} & $=$ & \includegraphics[scale=1,align=c]{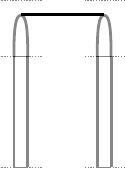} .

 \end{tabular}\end{center}

\end{enumerate}

\end{df}

An adjunction in a $3$-category is a pair of $1$-morphisms, together with unit and counit $2$-morphisms satisfying the snake relations up to isomorphism. By picking inverse pairs of $3$-morphisms - which are usually called \textbf{cusp $3$-morphisms} - implementing these snake relations and then considering the relations implied in witnessing that these are in fact inverse pairs of $3$-morphisms, one obtains all of the cells in the above presentation, except for the final one. This is the well known \textbf{swallowtail relation}. In \cite{verity_thesis} and \cite{gurski} the definitions of coherent adjunction include two swallowtail relations, but in \cite{piotr} the author shows that one follows from the other. We will give a string diagram proof of this fact.

Theorem \ref{main} will be a consequence of the following Lemma.

\begin{lmm}\label{pullback}

Given a $3$-category $\CC$, the map $$\Map(\Adj_{(3,1)},\CC)\to\Map(\Adj_{(3,1)},\h_2{\CC})\times_{\Map^L(\theta^{(1)}, \h_2{\CC})}\Map^L(\theta^{(1)}, \CC)$$ induced by the square $$\xymatrix{\Map(\Adj_{(3,1)},\CC)\ar[d]\ar[r]^-{E_l} & \Map^L(\theta^{(1)}, \CC)\ar[d] \\ \Map(\Adj_{(3,1)},\h_2{\CC})\ar[r]_-{E_l} & \Map^L(\theta^{(1)}, \h_2{\CC})}$$ is a weak equivalence.

\end{lmm}

First we explain how Theorem \ref{main} follows from this Lemma. We need one more standard definition and a another standard Lemma.

\begin{df}

A map of $n$-groupoids $f:X\to Y$ is a \textbf{trivial fibration} if it is both a fibration and a weak equivalence.

\end{df}

\begin{lmm}\label{triv_fib}

Consider a pullback

$$\xymatrix{X\times_B Y\ar[r]\ar[d] & Y \ar[d]^g \\ X\ar[r]_-{f} & B, }$$ where $f$ is a trivial fibration. Then $\pi_2:X\times_B Y \to Y$ is a trivial fibration.

\end{lmm}

\begin{proof}

The map $\pi_2$ is a fibration, because a lifting problem for $\pi_2$ reduces to a lifting problem for $f$, which is a fibration. The map $\pi_2$ is essentially surjective (on objects) because $f$ is essentially surjective. Finally, $\pi_2^{-1}(y)=f^{-1}(g(y))$ is weakly contractible, because $f$ is a fibration and a weak equivalence.
\end{proof}

\begin{proof}[Proof of Theorem \ref{main}]

Note that the bottom map in the square from Lemma \ref{pullback} is the same as $$\Map(\Adj_{(2,1)},\h_2{\CC})\to\Map^L(\theta^{(1)}, \h_2{\CC})$$ which is a fibration by Theorem \ref{fibration} and a weak equivalence of $2$-groupoids by Proposition \ref{2d}. Therefore, applying Lemma \ref{triv_fib}, the map $$\pi_2:\Map(\Adj_{(3,1)},\h_2{\CC})\times_{\Map^L(\theta^{(1)}, \h_2{\CC})}\Map^L(\theta^{(1)},\CC)\to\Map^L(\theta^{(1)},\CC)$$ is a trivial fibration of $3$-groupoids. By Lemma \ref{pullback}, the composite $\Map(\Adj_{(3,1)},\CC)\to\Map^L(\theta^{(1)},\CC)$ is a weak equivalence.
\end{proof}

So we need to prove Lemma \ref{pullback}, which is to say we need to show that the map $$\Map(\Adj_{(3,1)},\CC)\to\Map(\Adj_{(3,1)},\h_2{\CC})\times_{\Map^L(\theta^{(1)}, \h_2{\CC})}\Map^L(\theta^{(1)}, \CC)$$ is a weak equivalence of $3$-groupoids. Since this map is a fibration, it is enough to show that it is surjective on objects and has weakly contractible fibres.

\subsection{Surjective on objects}

We show that the map $$\Map(\Adj_{(3,1)},\CC)\to\Map(\Adj_{(3,1)},\h_2{\CC})\times_{\Map^L(\theta^{(1)}, \h_2{\CC})}\Map^L(\theta^{(1)}, \CC)$$ is surjective on objects.

\begin{df}

We define $\overline{\Adj}_{(3,1)}$ to be the presentation obtained from $\Adj_{(3,1)}$ by removing the swallowtail relation.

\end{df}

\begin{lmm}

Let $\CC$ be a $3$-category and consider a functor $$\overline{F}:\overline{\Adj}_{(3,1)}\to\CC.$$ Then there exists a functor $F:\Adj_{(3,1)}\to\CC$ which agrees with $\overline{F}$ on all cells in $\overline{\Adj}_{(3,1)}$ except for $C_l$ and $C_l^{-1}$.

\end{lmm}

\begin{proof}

We use the same notation for the image of a cell under $\overline{F}$ as for the cell itself. We define $$F(C_l)=\includegraphics[scale=1,align=c]{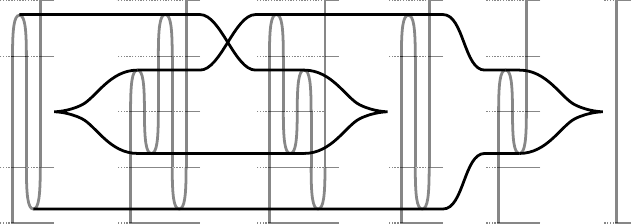}$$ and $$F(C_l^{-1})=\includegraphics[scale=1,align=c]{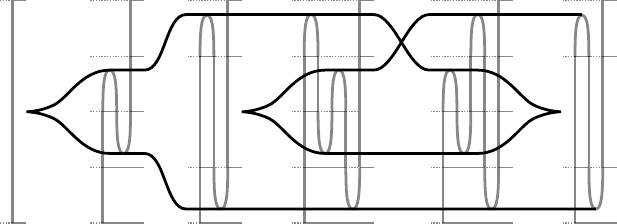}.$$ These are clearly inverse to each other and the following is a proof that the swallowtail relation holds with this choice of cusps:

\begin{center}\includegraphics[scale=1]{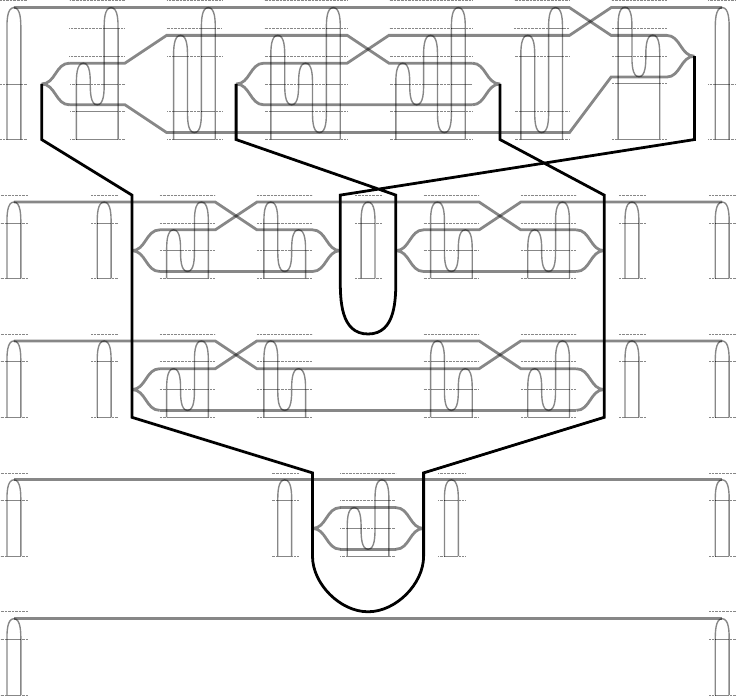} .\end{center}
\end{proof}

\begin{lmm}

The map $$\Map(\Adj_{(3,1)},\CC)\to\Map(\Adj_{(3,1)},\h_2{\CC})\times_{\Map^L(\theta^{(1)}, \h_2{\CC})}\Map^L(\theta^{(1)}, \CC)$$ is surjective on objects.

\end{lmm}

\begin{proof}

Given a map $F:\Adj_{(3,1)}\to \h_2{\CC}$ we must lift it to $\Adj_{(3,1)}\to\CC$. From $F:\Adj_{(3,1)}\to \h_2{\CC}$ we can build a map $\overline{F}:\overline{\Adj}_{(3,1)}\to\CC$ whose image in $\Map(\Adj_{(3,1)},\h_2{\CC})$ is $F$, so we can apply the previous Lemma.
\end{proof}

\subsection{Additional swallowtail relations}

Denote by $(\SW)$ the swallowtail relation in $\Adj_{(3,1)}$. We show that there are three additional swallowtail relations $(\overline{\SW})$, $(\SW_2)$ and $(\overline{\SW_2})$ which hold in $F(\Adj_{(3,1)})$. We will therefore use all four relations freely in the rest of the paper. The relation $(\SW_2)$ is usually included in the definition of a coherent adjunction. The fact that it follows from $(\SW)$ is originally due to \cite{piotr}, in the context of duals in monoidal bicategories. We present here a new string diagram proof. The sources of the relations $(\overline{\SW})$ and $(\overline{\SW_2})$ are inverse to those of $(\SW)$ and $(\SW_2)$, so they follow trivially.

\begin{df}

We define $\Adj^{+}_{(3,1)}$ to be the presentation obtained from $\Adj_{(3,1)}$ by adding the relations

\begin{center}
 \begin{center}\begin{tabular}{llllll}

$(\overline{\SW})$ & \includegraphics[scale=1,align=c]{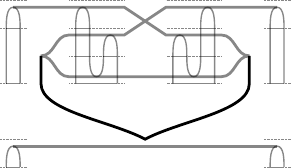} & $:$ & \includegraphics[scale=1,align=c]{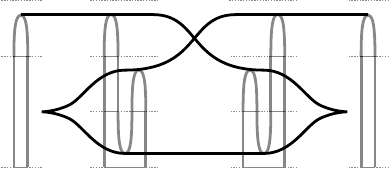} & $=$ & \includegraphics[scale=1,align=c]{adj/pres/id_u.pdf} ;

\\

\\

$(\SW_2)$ & \includegraphics[scale=1,align=c]{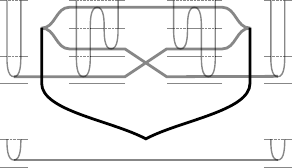} & $:$ & \includegraphics[scale=1,align=c]{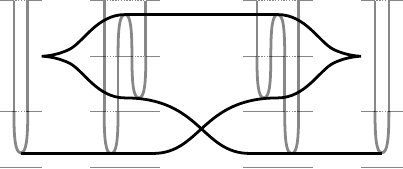} & $=$ & \includegraphics[scale=1,align=c]{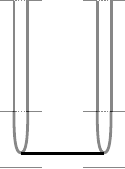} ;

\\

\\

$(\overline{\SW_2})$ & \includegraphics[scale=1,align=c]{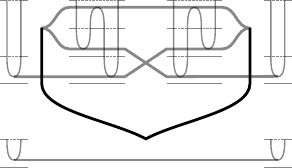} & $:$ & \includegraphics[scale=1,align=c]{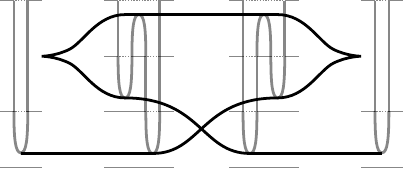} & $=$ & \includegraphics[scale=1,align=c]{adj/pres/id_c.pdf} .

 \end{tabular}\end{center}

 \end{center}

\end{df}

\begin{p}

We have $F(\Adj_{(3,1)}^{+})=F(\Adj_{(3,1)})$.

\end{p}

\begin{proof}

We have to show that the extra relations are already satisfied in $F(\Adj_{(3,1)})$. The source of $\overline{\SW}$ is inverse to the source of the swallowtail relation in $\Adj_{(3,1)}$, so it follows. We have the following proof for $\SW_2$:

\begin{center}

\includegraphics[scale=1]{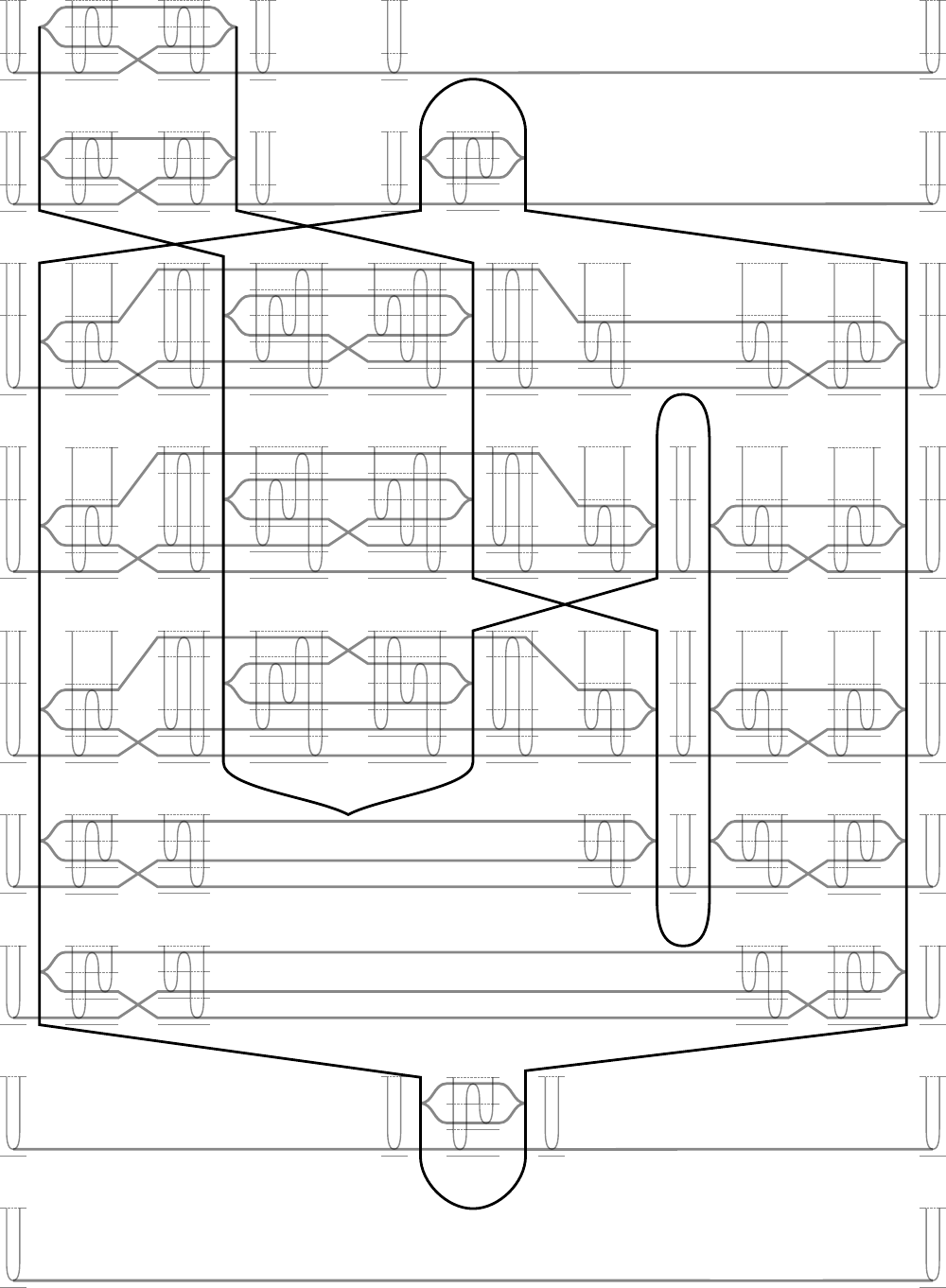} .

\end{center}

Finally, the source of $\overline{\SW_2}$ is inverse to the source of $\SW_2$.
\end{proof}

\subsection{Fibres are connected}

Now we prove that the fibres are connected. Given $F,G\in\Map(\Adj_{(3,1)}, \CC)$ with the same image under the map $$\Map(\Adj_{(3,1)}, \CC)\to\Map(\Adj_{(3,1)},\h_2{\CC})\times_{\Map^L(\theta^{(1)}, \h_2{\CC})}\Map^L(\theta^{(1)}, \CC),$$ we need to find an equivalence $\alpha:F\to G$ in $\Map(\Adj_{(3,1)},\CC)$ that maps to the identity. Equivalently, given $F,G\in\Map(\Adj_{(3,1)}, \CC)$ with $F=G$ on the $2$-skeleton of $\Adj_{(3,1)}$, we need to find an equivalence $\alpha:F\to G$ in $\Map(\Adj_{(3,1)},\CC)$ which is the identity on the $1$-skeleton of $\Adj_{(3,1)}$.

We use {\red} and {\blue} for the images of cells under functors $F$ and $G$, respectively. When $F$ and $G$ agree on a cell, we use black. We use {\green} for the values of $\alpha$.

\begin{lmm}\label{acldef}

Let $\CC$ be a $3$-category. Given $F,G\in\Map(\Adj_{(3,1)},\CC)$ with $F=G$ on $\sk_2(\Adj_{(3,1)})$ there exists an equivalence $\alpha:F\to G$ in $$\Map(\sk_2(\Adj_{(3,1)})\cup\{C_l\},\CC)$$ which is the identity on $\sk_1(\Adj_{(3,1)})\cup\{u\}$.

\end{lmm}

\begin{proof}

We need to construct $\alpha_c=\includegraphics[scale=1,align=c]{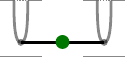}$ and show that the relation $\alpha_{C_l}$ is satisfied: \begin{center}\includegraphics[scale=1.5,align=c]{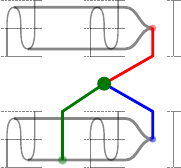} .\end{center} We define \begin{center}$\alpha_c:=$\includegraphics[scale=1,align=c]{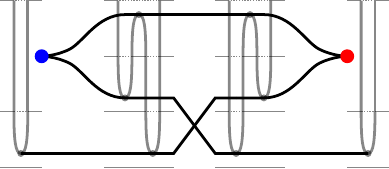}\end{center} and then we have the following proof of $\alpha_{C_l}:$ \begin{center}\includegraphics[scale=1]{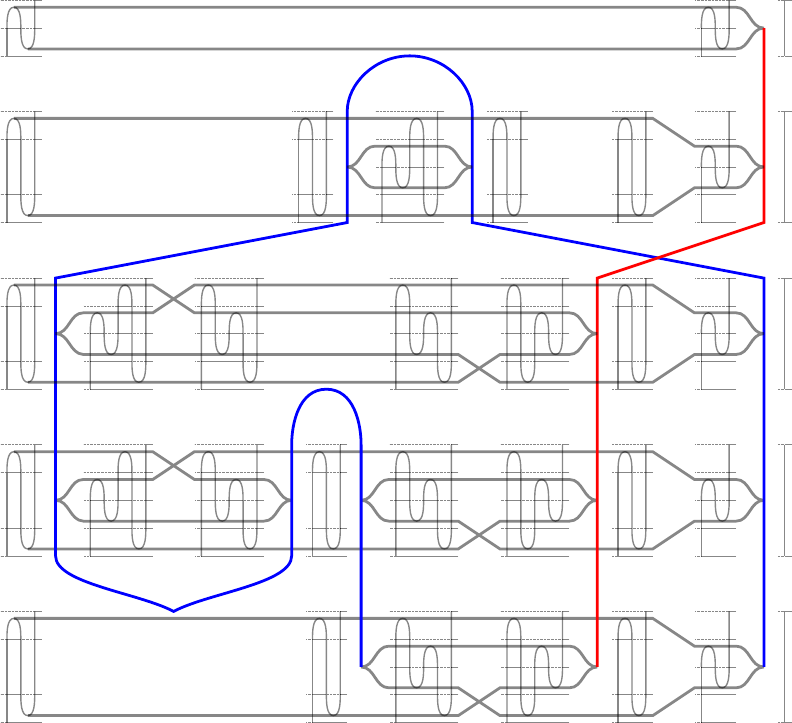} .\end{center}
\end{proof}

\begin{lmm}

Let $\CC$ be a $3$-category. Given $F,G\in\Map(\Adj_{(3,1)},\CC)$ with $F=G$ on $\sk_2(\Adj_{(3,1)})\cup\{C_l\}$, we have $F=G$.

\end{lmm}

\begin{proof}

First we have $F(C_l^{-1})=G(C_l^{-1})$ because both are both inverse to $F(C_l)=G(C_l)$. Then we have $F(C_r)=G(C_r)$ by \begin{center}\includegraphics[scale=1]{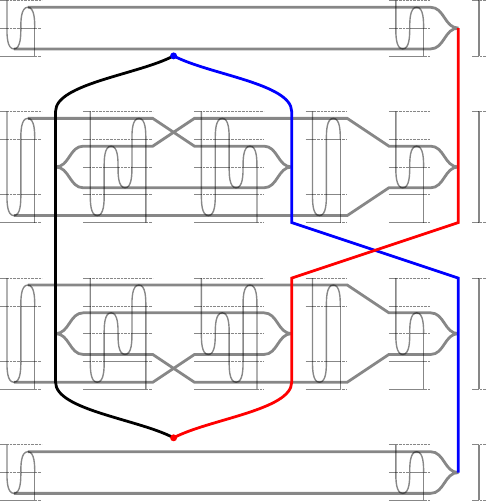} .\end{center} Finally, this implies $F(C_r^{-1})=G(C_r^{-1})$.
\end{proof}

\begin{lmm}

Let $\CC$ be a $3$-category. The fibres of $$\Map(\Adj_{(3,1)},\CC)\to\Map(\Adj_{(3,1)},\h_2{\CC})\times_{\Map^L(\theta^{(1)}, \h_2{\CC})}\Map^L(\theta^{(1)}, \CC)$$ are connected.

\end{lmm}

\begin{proof}

Given $F,G\in\Map(\Adj_{(3,1)}, \CC)$ with $F=G$ on $\{X,Y,l,r,u,c\}$, we need to find an equivalence $$\alpha:F\to G$$ in $\Map(\Adj_{(3,1)},\CC)$ which is the identity on $\{X,Y,l,r\}$. By Lemma \ref{acldef}, there exists an equivalence $$\alpha:F\to G$$ in $\Map(\{X,Y,l,r,u,c,C_l\},\CC)$, which is the identity on $\{X,Y,l,r,u\}$. We can lift this to an equivalence $F\to F_1$ in $\Map(\Adj_{(3,1)},\CC)$. Then $F_1=G$ on $\{X,Y,l,r,u,c,C_l\}$, so $F_1=G$.
\end{proof}

\subsection{Fibres are $1$-connected}

Now we prove that the fibres are $1$-connected. Given $F\in\Map(\Adj_{(3,1)},\CC)$ and an equivalence $\alpha:F\to F$ which is the identity on the $1$-skeleton, we need to find a $2$-equivalence $m:\alpha\to \Id$ in $\Map(\Adj_{(3,1)},\CC)$ which is the identity on $\{X,Y,l\}$. We use {\green} and {\orange} for the values of $\alpha$ and $m$, respectively.

\begin{lmm}\label{mudef}

Let $\CC$ be a $3$-category. Given a $1$-morphism $\alpha:F\to F$ in $\Map(\Adj_{(3,1)},\CC)$ such that $\alpha$ is the identity on $\sk_1(\Adj_{(3,1)})$, there exists a $2$-morphism $m:\alpha\to\Id_F$ in $\Map(\sk_1(\Adj_{(3,1)})\cup\{u\},\CC)$ which is the identity on $\{X,Y,l\}$.

\end{lmm}

\begin{proof}

Denote $\alpha_u:F(u)\to F(u)$ by \begin{center}\includegraphics[scale=1,align=c]{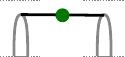}.\end{center} We want to define a $3$-morphism $m_r:\Id_r\to\Id_r$, which we denote by \includegraphics[scale=1]{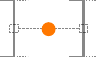}, such that the relation \begin{center}$m_u:$\includegraphics[scale=1,align=c]{adj/2morph/m_u_lhs_special.pdf}$=$ \includegraphics[scale=1,align=c]{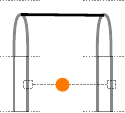}\end{center} is satisfied. Define \begin{center}$m_r:=$\includegraphics[scale=1,align=c]{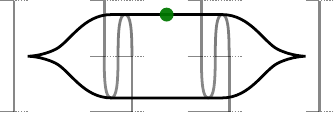}\end{center} and then the following is a proof of $m_u:$ \begin{center}\includegraphics[scale=1]{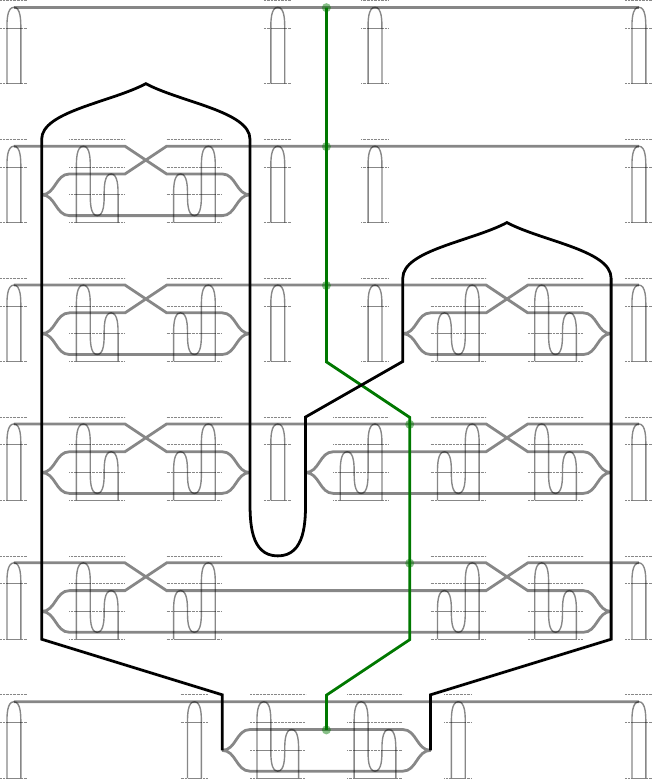} .\end{center}
\end{proof}

\begin{lmm}\label{mcdef}

Let $\CC$ be a $3$-category. Given a $1$-morphism $\alpha:F\to F$ in $\Map(\Adj_{(3,1)},\CC)$ such that $\alpha$ is the identity on $\sk_1(\Adj_{(3,1)})\cup\{u\}$, we have $\alpha=\Id$.

\end{lmm}

\begin{proof}

Denote $\alpha_c:F(c)\to F(c)$ by \begin{center}\includegraphics[scale=1,align=c]{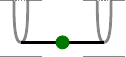}\end{center} and consider the relation $\alpha_{C_l^{-1}}:$\begin{center}\includegraphics[scale=1.5]{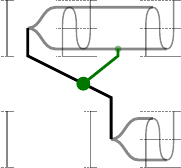} .\end{center} The following is a proof that $\alpha_c=\Id_c:$ \begin{center}\includegraphics[scale=1]{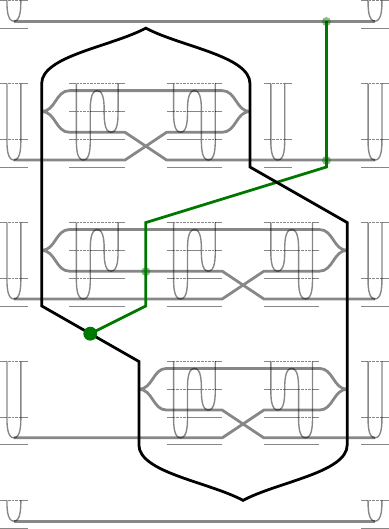} .\end{center}
\end{proof}

\begin{lmm}

Let $\CC$ be a $3$-category. The fibres of $$\Map(\Adj_{(3,1)},\CC)\to\Map(\Adj_{(3,1)},\h_2{\CC})\times_{\Map^L(\theta^{(1)}, \h_2{\CC})}\Map^L(\theta^{(1)}, \CC)$$ are $1$-connected.

\end{lmm}

\begin{proof}

Consider $F\in\Map(\Adj_{(3,1)},\CC)$ and an equivalence $\alpha:F\to F$ which is the identity on $\{X,Y,l,r\}$. By Lemma \ref{mudef}, there exists a $2$-morphism $m:\alpha\to\Id_F$ in $\Map(\{X,Y,l,r,u\},\CC)$ which is the identity on $\{X,Y,l\}$. We extend this to an equivalence $\alpha\to\alpha_1$ in $\Map(\Adj_{(3,1)},\CC)$ with $\alpha_1=\Id$ on $\{X,Y,l,r,u\}$. Then, by Lemma \ref{mcdef}, we have $\alpha_1=\Id$, so we have an equivalence $\alpha\to\Id$ which is the identity on $\{X,Y,l\}$.
\end{proof}

\subsection{Fibres are $2$-connected}

Now we prove that the fibres of $$\Map(\Adj_{(3,1)},\CC)\to\Map(\Adj_{(3,1)},\h_2{\CC})\times_{\Map^L(\theta^{(1)}, \h_2{\CC})}\Map^L(\theta^{(1)}, \CC)$$ are $2$-connected.

\begin{lmm}\label{aardef}

Let $\CC$ be a $3$-category. Given a $2$-morphism $m:\Id_F\to \Id_F$ in $\Map(\Adj_{(3,1)},\CC)$ such that $m$ is the identity on $\{X,Y,l\}$, we have $m=\Id$.

\end{lmm}

\begin{proof}

Denote \begin{center}$m_r:=$ \includegraphics[scale=1,align=c]{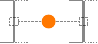}\end{center} and consider the relation $m_c:$ \begin{center} \includegraphics[scale=1.5]{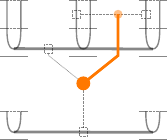} .\end{center} We have the following proof that $m_r= \Id^{(2)}_{r}:$ \begin{center}\includegraphics[scale=1]{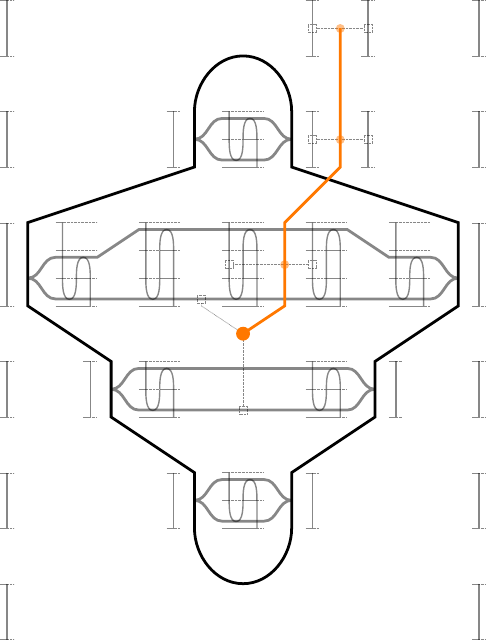} .\end{center}
\end{proof}

\begin{lmm}

Let $\CC$ be a $3$-category. The fibres of $$\Map(\Adj_{(3,1)},\CC)\to\Map(\Adj_{(3,1)},\h_2{\CC})\times_{\Map^L(\theta^{(1)}, \h_2{\CC})}\Map^L(\theta^{(1)}, \CC)$$ are $2$-connected.

\end{lmm}

\begin{proof}

This follows from the previous Lemma.
\end{proof}

\begin{lmm}

Let $\CC$ be a $3$-category. The fibres of $$\Map(\Adj_{(3,1)},\CC)\to\Map(\Adj_{(3,1)},\h_2{\CC})\times_{\Map^L(\theta^{(1)}, \h_2{\CC})}\Map^L(\theta^{(1)}, \CC)$$ are weakly contractible.

\end{lmm}

\begin{proof}

Since morphisms in the fibres have to restrict to the identity on the $0$-skeleton of $\Adj_{(3,1)}$, the fibres are $2$-groupoids. Since we have shown they are $2$-connected, they are therefore weakly contractible.
\end{proof}

This concludes the proof of Lemma \ref{pullback} and therefore of Theorem \ref{main}.

\subsection*{Acknowledgements}

I am grateful to Pedro Boavida, Christopher Douglas, John Huerta and Roger Picken, as well as the two anonymous reviewers, for many useful comments and suggestions. This work was partially supported by the FCT project grant \textbf{Higher Structures and Applications}, PTDC/MAT-PUR/31089/2017. The work that lead to this paper was also carried out while the author was visiting the \textbf{Max Planck Institute for Mathematics} and later a postdoc with the RTG 1670 \textbf{Mathematics inspired by String Theory and Quantum Field Theory} at the University of Hamburg.

\bibliographystyle{plainnat}
\bibliography{refs}

\begin{thebibliography}{21}
\providecommand{\natexlab}[1]{#1}
\providecommand{\url}[1]{\texttt{#1}}
\expandafter\ifx\csname urlstyle\endcsname\relax
  \providecommand{\doi}[1]{doi: #1}\else
  \providecommand{\doi}{doi: \begingroup \urlstyle{rm}\Url}\fi

\bibitem[Ara and Lucas(2020)]{al20}
Dimitri Ara and Maxime Lucas.
\newblock {The folk model category structure on strict $\omega$-categories is
  monoidal}.
\newblock \emph{{Theory and Applications of Categories}}, 35\penalty0
  (21):\penalty0 745--808, 2020.
\newblock URL \url{http://www.tac.mta.ca/tac/volumes/35/21/35-21abs.html}.

\bibitem[Ara and M{\'e}tayer(2011)]{am11}
Dimitri Ara and Fran{\c c}ois M{\'e}tayer.
\newblock {The Brown-Golasinski model structure on strict $\infty$-groupoids
  revisited}.
\newblock \emph{{Homology, Homotopy and Applications}}, 13\penalty0
  (1):\penalty0 121--142, 2011.
\newblock \doi{10.4310/HHA.2011.v13.n1.a6}.

\bibitem[Ara\'{u}jo(2017)]{thesis}
Manuel Ara\'{u}jo.
\newblock \emph{Coherence for 3-dualizable objects}.
\newblock PhD thesis, University of Oxford, 2017.
\newblock URL
  \url{https://ora.ox.ac.uk/objects/uuid:a4b8f8de-a8e3-48c3-a742-82316a7bd8eb}.

\bibitem[Ara\'{u}jo(2020)]{adjoints1}
Manuel Ara\'{u}jo.
\newblock String diagrams for 4-categories and fibrations of mapping
  4-groupoids, 2020.
\newblock URL \url{https://arxiv.org/abs/2012.03797}.

\bibitem[Ara\'{u}jo(2022)]{sesqui}
Manuel Ara\'{u}jo.
\newblock Simple string diagrams and $n$-sesquicategories, 2022.
\newblock URL \url{https://arxiv.org/abs/2202.09293}.

\bibitem[Bar and Vicary(2017)]{Jamie}
Krzysztof Bar and Jamie Vicary.
\newblock Data structures for quasistrict higher categories.
\newblock In \emph{Proceedings of 32nd Annual ACM/IEEE Symposium on Logic in
  Computer Science (LICS 2017)}. IEEE Computer Society, 2017.
\newblock \doi{10.1109/LICS.2017.8005147}.

\bibitem[Bar et~al.(2016)Bar, Kissinger, and Vicary]{globular}
Krzysztof Bar, Aleks Kissinger, and Jamie Vicary.
\newblock {Globular: An Online Proof Assistant for Higher-Dimensional
  Rewriting}.
\newblock In \emph{1st International Conference on Formal Structures for
  Computation and Deduction (FSCD 2016)}, volume~52 of \emph{Leibniz
  International Proceedings in Informatics (LIPIcs)}, pages 34:1--34:11, 2016.
\newblock \doi{10.4230/LIPIcs.FSCD.2016.34}.

\bibitem[Barrett et~al.(2012)Barrett, Meusburger, and Schaumann]{surface}
John~W. Barrett, Catherine Meusburger, and Gregor Schaumann.
\newblock Gray categories with duals and their diagrams, 2012.
\newblock URL \url{https://arxiv.org/abs/1211.0529}.

\bibitem[Bartlett(2009)]{bruce}
Bruce Bartlett.
\newblock \emph{On unitary 2-representations of finite groups and topological
  quantum field theory}.
\newblock PhD thesis, University of Sheffield, 2009.
\newblock URL \url{https://arxiv.org/abs/0901.3975}.

\bibitem[Brown and Golasinski(1989)]{bg89}
Ronald Brown and Marek Golasinski.
\newblock A model structure for the homotopy theory of crossed complexes.
\newblock \emph{Cahiers de Topologie et G\'{e}om\'{e}trie Diff\'{e}rentielle
  Cat\'{e}goriques}, 30\penalty0 (1):\penalty0 61--82, 1989.
\newblock URL \url{http://eudml.org/doc/91432}.

\bibitem[Crans(1995)]{crans95}
Sjoerd~E. Crans.
\newblock Pasting schemes for the monoidal biclosed structure on $\omega$-cat,
  1995.
\newblock URL
  \url{http://citeseerx.ist.psu.edu/viewdoc/summary?doi=10.1.1.56.8738}.

\bibitem[Gurski(2012)]{gurski}
Nick Gurski.
\newblock Biequivalences in tricategories.
\newblock \emph{Theory and Applications of Categories}, 26\penalty0
  (14):\penalty0 349--384, 2012.
\newblock URL \url{http://www.tac.mta.ca/tac/volumes/26/14/26-14abs.html}.

\bibitem[Johnson-Freyd and Scheimbauer(2017)]{freyd_scheim}
Theo Johnson-Freyd and Claudia Scheimbauer.
\newblock (op)lax natural transformations, twisted quantum field theories, and
  {"}even higher{"} morita categories.
\newblock \emph{Advances in Mathematics}, 307:\penalty0 147 -- 223, 2017.
\newblock \doi{10.1016/j.aim.2016.11.014}.

\bibitem[Lafont et~al.(2010)Lafont, M{\'e}tayer, and Worytkiewicz]{folk}
Yves Lafont, Fran{\c c}ois M{\'e}tayer, and Krzysztof Worytkiewicz.
\newblock {A folk model structure on omega-cat}.
\newblock \emph{{Advances in Mathematics}}, 224\penalty0 (3):\penalty0
  1183--1231, 2010.
\newblock \doi{10.1016/j.aim.2010.01.007}.

\bibitem[Leinster(2004)]{OperadsCats}
Tom Leinster.
\newblock \emph{Higher Operads, Higher Categories}.
\newblock London Mathematical Society Lecture Note Series. Cambridge University
  Press, 2004.
\newblock \doi{10.1017/CBO9780511525896}.

\bibitem[Lurie(2008)]{lurie}
Jacob Lurie.
\newblock On the classification of topological field theories.
\newblock \emph{Current Developments in Mathematics}, 2008:\penalty0 129 --
  280, 2008.
\newblock \doi{10.4310/CDM.2008.v2008.n1.a3}.

\bibitem[Pstr\k{a}gowski(2022)]{piotr}
Piotr Pstr\k{a}gowski.
\newblock On dualizable objects in monoidal bicategories.
\newblock \emph{Theory and Applications of Categories}, 38\penalty0
  (9):\penalty0 257--310, 2022.
\newblock URL \url{http://www.tac.mta.ca/tac/volumes/38/9/38-09abs.html}.

\bibitem[Riehl and Verity(2016)]{riehlverity}
Emily Riehl and Dominic Verity.
\newblock Homotopy coherent adjunctions and the formal theory of monads.
\newblock \emph{Advances in Mathematics}, 286:\penalty0 802--888, 2016.
\newblock \doi{10.1016/j.aim.2015.09.011}.

\bibitem[Schanuel and Street(1986)]{schstr}
Stephen Schanuel and Ross Street.
\newblock The free adjunction.
\newblock \emph{Cahiers de Topologie et G\'eom\'etrie Diff\'erentielle
  Cat\'egoriques}, 27\penalty0 (1):\penalty0 81--83, 1986.
\newblock URL \url{http://www.numdam.org/item/CTGDC_1986__27_1_81_0/}.

\bibitem[Schommer-Pries(2009)]{CSPThesis}
Christopher Schommer-Pries.
\newblock \emph{The Classification of Two-Dimensional Extended Topological
  Field Theories}.
\newblock PhD thesis, University of California, Berkeley, 2009.
\newblock URL \url{https://arxiv.org/abs/1112.1000}.

\bibitem[Verity(1992)]{verity_thesis}
Dominic Verity.
\newblock \emph{Enriched categories, internal categories and change of base}.
\newblock PhD thesis, University of Cambridge, 1992.
\newblock URL
  \url{http://www.tac.mta.ca/tac/reprints/articles/20/tr20abs.html}.

\end{thebibliography}

\end{document}